\documentclass[12pt]{amsart}

\usepackage{amscd}

\usepackage{comments}
\newComments\DL{DL}{red}
\newComments\BJ{BJ}{blue}
\newComments\IS{IS}{red}

\usepackage{soul}
\usepackage{lscape}

\usepackage{DLdef1}

\newcommand {\size}{{\text{size}}}
\newcommand {\foo}{{\mathfrak{oo}}}

\newcommand {\foc}{{\mathfrak{oc}}}

\newcommand {\fbgl}{{\mathfrak{bgl}}}

\newcommand {\fwk}{{\mathfrak{wk}}}
\newcommand{\uN}{\underline{N}}
\newcommand{\un}{\underline{N}}
\newcommand{\del}{\partial}

\begin{document}

\title[New simple Lie algebras in characteristic $2$]
{Simple prolongs of the non-positive parts of graded Lie algebras
with Cartan matrix in characteristic $2$}

\author{Sofiane Bouarroudj${}^1$, Pavel Grozman${}^2$, Alexei
Lebedev${}^2$, Dimitry Leites${}^3$,\\ Irina Shchepochkina${}^4$}

\address{${}^1$New York University Abu Dhabi,
Division of Science and Mathematics, P.O. Box 129188, United Arab
Emirates; sofiane.bouarroudj@nyu.edu}\address{${}^2$Equa Simulation
AB, Stockholm, Sweden;
pavel.grozman@bredband.net}\address{${}^3$Department of Mathematics,
Stockholm University, Ros\-lagsv. 101, Kr\"aft\-riket hus 6, SE-106
91 Stockholm, Sweden; mleites@math.su.se}\address{${}^4$Independent
University of Moscow, Bolshoj Vlasievsky per, dom 11, RU-119 002
Moscow, Russia; irina@mccme.ru}

\keywords {Cartan prolongation, Lie superalgebra, characteristic 2,
Kostrikin-Shafarevich conjecture}

\subjclass{17B50, 70F25}

\begin{abstract} Over an algebraically closed fields, an alternative
to the method due to Kostrikin and Shafarevich was recently
suggested. It produces all known simple finite dimensional Lie
algebras in characteristic $p>2$. For $p=2$, we investigate one of
the steps of this method, interpret several other simple Lie
algebras, previously known only as sums of their components, as Lie
algebras of vector fields. One new series of exceptional simple Lie
algebras is discovered, together with its ``hidden supersymmetries".

In characteristic 2, certain simple Lie algebras are
``desuperizations" of simple Lie superalgebras. Several simple Lie
algebras we describe as results of generalized Cartan prolongation
of the non-positive parts, relative a simplest (by declaring degree
of just one pair of root vectors corresponding to opposite simple
roots nonzero) grading by integers, of Lie algebras with Cartan
matrix are ``desuperizations" of characteristic 2 versions of
complex simple exceptional vectorial Lie superalgebras. We list the
Lie superalgebras (some of them new) obtained from the Lie algebras
considered by declaring certain generators odd.

One of the simple Lie algebras obtained is the prolong relative to a
non-simplest grading, so the classification to be obtained might be
more involved than we previously thought.
\end{abstract}

\thanks{SB is partly supported by grant AD 65 NYUAD}


\maketitle

\markboth{\itshape Sofiane Bouarroudj\textup{,} Pavel
Grozman\textup{,} Alexei Lebedev\textup{,} Dimitry Leites\textup{,}
Irina Shchepochkina}{{\itshape New simple Lie algebras in
characteristic $2$}}

\thispagestyle{empty}

\section{Introduction}
Hereafter, $\Kee$ is an algebraically closed field of characteristic
$2$ unless indicated other\-wise. For background, see
\cite{BGL1,LeP} and \S2 which has less examples but a few more
clarifications.

\ssec{Main results of this paper} For $p=2$ and the non-positive
part of each Lie algebra $\fg(A)$ with indecomposable Cartan matrix
$A$ (classified in \cite{WK, BGL1}) for the simplest $\Zee$-gradings
of $\fg(A)$ and $A$ of size\footnote{This limitation is imposed by
in-built weaknesses of \textit{Mathematica} on which the package
\textit{SuperLie} we use for computer-aided studies is based. To
advance without computer aid is hopeless. However, we conjecture
that there are no new simple Lie algebras obtained as CTS prolongs
of the non-positive parts of algebras $\fg(A)$ for Cartan matrices
of larger size.} $\leq 4$, we compute the CTS prolong $(\fg_-,
\fg_0)_{*,\un}$. If $\fg(A)$ is not simple, we consider prolongs of
the non-positive part of both $\fg(A)/\fc$, where $\fc$ is the
center, and the simple derived of $\fg(A)/\fc$. We denote by $\bf F$
the desuperization functor, the one that forgets parity of the Lie
superalgebra turning it into a Lie algebra (recall that $p=2$).
Considering the CTS-prolongs of the non-positive or negative parts
(relative certain particular $\Zee$-gradings) of exceptional
(discovered by Weisfeiler and Kac) Lie algebras $\fwk(3; a)/\fc$,
$\fwk'(3;a)/\fc$ and $\fwk(4; a)$, we obtain several simple Lie
algebras as \emph{desuperizations} of certain Lie
superalgebras\footnote{Brown did not give any interpretation of the
three series of Lie algebras he described only in components. For
$\fwk'(3; a)/\fc$, Brown \cite{Bro} was the first to consider
prolongs of its non-positive part for one of gradings corresponding
to one of several Cartan matrices of $\fwk(3; a)$. We interpret
another Brown's series, $D_4(3;\un)$, as a desuperization of the
exceptional simple vectorial Lie superalgebra $\fvle(3;\un|8)$
described for constrained values of the shearing vector $\un$. We
interpret the third Brown's series in \cite{BGLLS}.}.

These simple Lie algebras are new, but not to us: they are
\emph{desuperizations} of the characteristic 2 analogs of certain
simple exceptional complex vectorial Lie superalgebras, cf.
\cite{BGLS}. \emph{Completely new} are our \emph{interpretations} of
these desuperizations:

(a) ${\bf F}(\fm\fb(3;\un|8))$ as independent of parameter $a$
prolong of $\fwk(4; a)$,

(b) the main deform of the anti-bracket or Buttin superalgebra as
the prolong of $\fwk'(3;a)/\fc$ for a certain Cartan matrix and
grading of $\fwk(3;a)$;

(c) ${\bf F}(\fvle(3;\un|8))$ which under certain restrictions on
$\un$ turns into the Brown algebra $D_4(3;\un)$ being more general
otherwise.

Completely new are the exceptional simple Lie algebra
$\fii\fr(9;\un)$ and its superizations $\fii\fr(3;\un|6)$ and
$\fii\fr(5;\un|4)$; for $\un=(1,\dots, 1)$, they turn into
$\fo_\Pi'(8)/\fc$, $\fo_\Pi'(2|6)/\fc$ and $\fo_\Pi'(4|4)/\fc$,
respectively.

\emph{Hidden supersymmetries} of the Lie algebras $(\fg_-,
\fg_0)_{*,\un}$ are, by definition, the Lie superalgebras one can
obtain from $(\fg_-, \fg_0)_{*,\un}$ by declaring some of the
generators odd. Our description of the prolongs makes the
description of ``hidden supersymmetries" of the prolongs obvious and
explicit.

\ssec{The KSh method} Over the algebraically closed fields $\Kee$ of
characteristic $p\geq 7$, the Kostrikin-Shafarevich procedure for
obtaining all simple finite dimensional Lie algebras consists of the
following steps:
\begin{equation}\label{Kan}
\begin{minipage}[l]{14.5cm}
\hskip4 mm  1) for input, take the two types of simple complex Lie
algebras:

a) those of the form $\fg(A)$ for a Cartan matrix $A$,

b)  infinite dimensional vectorial Lie algebras with polynomial
coefficients;

\hskip 4 mm  2) among bases allowing integer structure constants
select certain ones with the ``smallest" constants (Chevalley bases
for the algebras of the form $\fg(A)$ and divided powers for
vectorial Lie algebras) thus getting particular $\Zee$-forms of
these complex Lie algebras;

\hskip 4 mm  3) tensor the $\Zee$-forms obtained at step 2) by
$\Kee$ over $\Zee$;

\hskip 4 mm  4) select a simple (and finite dimensional in the
vectorial Lie algebra case) subquotient, called a simple ``relative"
in what follows;

\hskip 4 mm  5) deform the results obtained at step 4);

\hskip 4 mm  6) classify isomorphisms between Lie algebras obtained
at earlier steps.
\end{minipage}
\end{equation}

In \cite{L}, conjectures describing ways to obtain all simple finite
dimensional Lie algebras and superalgebras over an algebraically
closed field $\Kee$ of characteristic 2 were offered. These
conjectures were sharpened lately, but not sufficiently. In this
paper, a sequel to \cite{BGL1, BGLS, BGL2, LeP}, we perform a step
towards these classifications along one of the ways indicated in
\cite{L, GL}; this helps to make the conjecture more precise.

\ssec{A reformulation of the procedures leading to the list
conjectured by Kostrikin and Shafarevich} Let
$\fg=\mathop{\oplus}\limits_{i\in \Zee}\fg_i$,
$\fg_-=\mathop{\oplus}\limits_{i\leq 0}\fg_i$ and
$\fg_0\subset\fder_0\fg_-$ a subalgebra preserving $\Zee$-grading of
$\fg_-$. For $p=5$, the Melikyan algebras are obtained by means of a
generalized Cartan prolongation of another type of pairs $(\fg_{-},
\fg_0)$  as compared with the pairs of the input for a generalized
Cartan prolongations listed in \eqref{Kan}. Actually, Melikyan's
examples, especially their interpretation as generalized prolongs of
the non-positive part of\footnote{We denote the exceptional algebras
$\fg(2)$, $\ff(4)$, etc. by analogy with $\fgl(n)$, and in order not
to confuse with the 2nd, 4th, $i$th. component $\fg_i$ of the
$\Zee$-graded algebra $\fg=\oplus \fg_i$.} $\fg(2)$, and the
Yamaguchi theorem \cite{Y} (we will recall it in \S2, see also a
more accessible than \cite{Y} paper \cite{GL}), hint at another
approach to the construction of simple modular Lie algebras, the one
we recall and apply in what follows,\footnote{Actually, for $p=2$,
the situation is even more involved and resembles that of simple
vectorial Lie superalgebras over $\Cee$, see \cite{LSh}. We consider
the cases with the input Lie algebras distinct from those considered
in this paper in \cite{BGLS, BGLLS}.} see also the steps listed in
\cite{L}.

The proof of the generalized KSh conjecture for $p>3$, mainly due to
Premet and Strade, is based on earlier results by Block and Wilson
(restricted case for $p>5$), and several other results by other
authors; for a final touch and history, see \cite{BGP}.

In \cite{GL}, we offered a procedure not only  leading to the same
list of simple modular Lie algebras for $p>3$ as produced by the KSh
procedure and Melikyan's examples, but which for $p=3$, clarified
several previously inexplicable examples, corrected earlier
findings, and produced several new simple Lie algebras. Namely, our
main ingredients are Lie algebras of the form $\fg(A)$ only, while
the main procedure is a generalization of Cartan prolongation
procedure, either complete or --- this is important! --- {\bf
partial}\footnote{For the definition of partial prolongations,
algorithm including,  see \cite{Shch}.}:
\begin{equation}\label{answ}
\begin{minipage}[l]{14.5cm}
for $p\geq 3$, up to deformations and passage to the derived
algebras and their quotients modulo center, the simple finite
dimensional modular Lie algebras are the results of the
Cartan-Tanaka-Shchepochkina (CTS) prolongations (either complete or
partial) of the non-positive parts (relative certain
$\Zee$-gradings) of the Lie algebras of the form $\fg(A)$ or their
derived.
\end{minipage}
\end{equation}

Thus, instead of the two types of  Lie algebras required by the KSh
procedure as the input, we need only one type of Lie algebras (of
the form $\fg(A)$) subjected to one type of constructions (CTS
prolongations). Subsequent passage to the derived (first or second),
factorization modulo center, and deforming are common features of
both approaches; however, selection of isomorphisms, especially
among the deforms, although common to both approaches, becomes much
more involved for $p=3$ and, especially, $p=2$, cf. \cite{KuCh,BLW,
BLLS}.

\sssec{The list of simple modular Lie algebras related to those of
the form $\fg(A)$ and their ``hidden supersymmetries"} For the
classification of finite dimensional Lie algebras $\fg(A)$ with
indecomposable Cartan matrix $A$ over algebraically closed fields
$\Kee$ of characteristic $p>0$, see \cite{WK} with corrections in
\cite{Sk1} and clarifications in \cite{BGL1}. In \cite{BGL1}, we
gave precise definitions of Cartan matrix and related notions
(Dynkin diagrams, Chevalley generators, and more) specific to the
super and characteristic $p>0$ cases, and classified finite
dimensional Lie superalgebras $\fg(A)$ with indecomposable Cartan
matrix $A$ over algebraically closed fields $\Kee$ of characteristic
$p>0$. Each finite dimensional Lie superalgebra $\fg(A)$ with an
indecomposable $A$ is either simple itself or $\fg^{(i)}/\fc$, where
$\fc$ is the center and $i=1$ or 2, is simple.

The answer in the case $p=2$ turned out to be very interesting:
\emph{the Lie algebras of the form $\fg(A)$ possess a hidden
supersymmetry}. More precisely,
\begin{equation}\label{hidsusy}
\begin{minipage}[l]{14cm}
Each finite dimensional Lie superalgebra $\fg(A)$ with
indecomposable Cartan matrix $A$ can be obtained from the Lie
algebra of the form $\fg(A)$ with the same $A$ by declaring any
number of pairs (positive and respective negative) of its Chevalley
generators odd.
\end{minipage}
\end{equation}

The vectorial Lie algebras possess same property (hidden
supersymmetry). Moreover,

(a) in order to understand what are all the analogs of orthogonal
and symplectic Lie algebras for $p=2$ (being interested in their
prolongs such as Lie algebras of Hamiltonian or contact vector
fields), we have to take into account super versions of orthogonal
and symplectic Lie algebras, namely, the periplectic Lie
superalgebras;

(b) several of seemingly new examples we obtained as CTS prolongs
are \emph{desuperizations} of certain characteristic 2 analogs of
exceptional simple complex vectorial Lie superalgebras.

Therefore, we have to give a necessary background concerning not
only Lie algebras but Lie superalgebras as well.

\sssec{If $p=2$, other inputs are needed for the CTS-procedure in
\eqref{answ}} For $p=2$, the procedure conjecturally leading to the
complete description of simple Lie algebras becomes much more
complicated than \eqref{answ}, see \cite{L}. In addition to the step
a) in the following list of steps \eqref{descr} leading to all
simple finite dimensional Lie algebras for $p>3$ and ---
conjecturally --- for $p=3$, we need at least the ingredients listed
in other steps.
\begin{equation}\label{descr}
\begin{minipage}[l]{14cm}\small

If $p=2$, in steps a)--d) we should consider {\bf not only simplest} $\Zee$-gradings:

\hskip 4 mm A) prolongs (complete and partial) of

\hskip 4 mm Aa) the non-positive parts of the Lie algebras of the
form $\fg(A)$ or their derived, where $A$ is indecomposable;

\hskip 4 mm Ab) the non-positive parts of the orthogonal Lie algebra
without Cartan matrix (or its first or second derived, or a central
extension thereof), see \cite{LeP};

\hskip 4 mm Ac) the non-positive parts of the Shen algebra and of
certain exceptional pairs $(\fg_{-1}, \fg_0)$, where  $\fg_{-1}$ is
a $\fg_0$-module;

\hskip 4 mm B) the results of application of the functor forgetting
superstructure to the $p=2$ analogs of Shchepochkina's simple
exceptional Lie superalgebras (partly listed in \cite{BGLS});

\hskip 4 mm C) deforming the results obtained at step A) and B);

\hskip 4 mm  D) classification of isomorphisms between Lie algebras
obtained at earlier steps (this  becomes even more involved for
$p=2$, cf.  \cite{BLW, KuCh}).
\end{minipage}
\end{equation}

\ssec{Related open problems and conjectures} In this paper we tackle
step a) of the main conjecture \eqref{descr} for the Cartan matrices
of size $\leq 4$ and for simplest $\Zee$-gradings only (with one
exception). We conjecture that prolongs of the non-positive parts of
the Lie algebras with Cartan matrices of larger size and more
complicated gradings (bar the above exception) return the initial
algebra (as in the generic cases of the Yamaguchi's theorem). To
investigate this conjecture is an important problem.

\section{Notation and the
background}

\ssec{What Lie superalgebra in characteristic $2$
is}\label{Ssalgin2}

Let us give a naive definition of a Lie superalgebra for $p=2$. (For
a scientific one, as a Lie algebra in the category of
supervarieties, needed, for example, for a rigorous study and
interpretation of odd parameters of deformations, see \cite{LSh}.)
We define a Lie superalgebra as a superspace
$\fg=\fg_\ev\oplus\fg_\od$ such that the even part $\fg_\ev$ is a
Lie algebra, the odd part $\fg_\od$ is a $\fg_\ev$-module (made into
the two-sided one by symmetry; more exactly, by {\it anti}-symmetry,
but if $p=2$, it is the same) and on $\fg_\od$ a {\it squaring}
(roughly speaking, the halved bracket) is defined as a map
\begin{equation}\label{squaring}
\begin{array}{c}
x\mapsto x^2\quad \text{such that $(ax)^2=a^2x^2$ for any $x\in
\fg_\od$ and $a\in \Kee$, and}\\
{}(x+y)^2-x^2-y^2\text{~is a bilinear form on $\fg_\od$ with values
in $\fg_\ev$.}
\end{array}
\end{equation}
(We use a minus sign, so the definition also works for $p\neq 2$.)
The origin of this operation is as follows: If $\Char \Kee\neq 2$,
then for any Lie superalgebra $\fg$ and any odd element
$x\in\fg_\od$, the Lie superalgebra $\fg$ contains the element $x^2$
which is equal to the even element $\frac12 [x,x]\in\fg_\ev$. It is
desirable to keep this operation for the case of $p=2$, but, since
it can not be defined in the same way, we define it separately, and
then define the bracket of odd elements to be (this equation is
valid for $p\neq 2$ as well):
\begin{equation}\label{bracket}
{}[x, y]:=(x+y)^2-x^2-y^2.
\end{equation}
We also assume, as usual, that
\begin{itemize}
  \item if $x,y\in\fg_\ev$, then $[x,y]$ is the bracket on the Lie algebra;
  \item if $x\in\fg_\ev$ and $y\in\fg_\od$, then
$[x,y]:=l_x(y)=-[y,x]=-r_x(y)$, where $l$ and $r$ are the left and
right $\fg_\ev$-actions on $\fg_\od$, respectively.
\end{itemize}

The Jacobi identity involving odd elements now takes the following
form:
\begin{equation}\label{JI}
~[x^2,y]=[x,[x,y]]\text{~for any~} x\in\fg_\od, y\in\fg.
\end{equation}
If $\Kee\neq \Zee/2\Zee$, we can replace the condition (\ref{JI}) on
two odd elements by a simpler one:
\begin{equation}\label{JIbest}
[x,x^2]=0\ \text{ for any $x\in\fg_\od$.}
\end{equation}

Because of the squaring, the definition of derived algebras should
be modified. For any Lie superalgebra $\fg$, set $\fg^{(0)}:=\fg$
and
\begin{equation}\label{deralg}
\fg':=[\fg,\fg]+\Span\{g^2\mid g\in\fg_\od\} ,\quad
\fg^{(i+1)}:=[\fg^{(i)},\fg^{(i)}]+\Span\{g^2\mid
g\in\fg^{(i)}_\od\}.
\end{equation}


An even linear map $r\colon  \fg\tto\fgl(V)$ is said to be a {\it
representation of the Lie superalgebra}\index{representation of the
Lie superalgebra} $\fg$ (and the superspace $V$ is said to be a {\it
$\fg$-module})\index{$\fg$-module} if
\begin{equation}\label{repres}
\begin{array}{l}
r([x, y])=[r(x), r(y)]\quad \text{ for any $x, y\in
\fg$;}\\
r(x^2)=(r(x))^2\text{~for any $x\in\fg_\od$.}
\end{array}
\end{equation}

\sssec{Examples: Lie superalgebras preserving non-degenerate
(anti-)sym\-met\-ric forms} We say that two bilinear forms $B$ and
$B'$ on a superspace $V$ are {\it equivalent} if there is an even
invertible linear map $M\colon V\tto V$ such that
\begin{equation}\label{eqform}
B'(x,y)=B(Mx,My) \text{~~for any~~}x,y\in V.
\end{equation}
{\bf We fix some basis in $V$ and identify a given bilinear form
with its Gram matrix in this basis; we also identify any linear
operator on $V$ with its supermatrix in a fixed basis}.

Then two bilinear forms (rather supermatrices) are equivalent if and
only if there is an even invertible matrix $M$ such that
\begin{equation}\label{eqformM}
B'=MBM^T,\text{~~where $T$ is for transposition}.
\end{equation}

A bilinear form $B$ on $V$ is said to be {\it
symmetric}\index{bilinear form! symmetric} if $B(v, w)=B(w, v)$ for
any $v, w\in V$; a bilinear form is said to be {\it anti-symmetric}
if $B(v, v)=0$ for any $v\in V$.\index{bilinear form!
anti-symmetric}

A homogeneous\footnote{Hereafter, as always in Linear Algebra in
superspaces, all formulas of linear algebra defined on homogeneous
elements only are supposed to be extended to arbitrary ones by
linearity.} linear map $F$ is said to preserve a bilinear form $B$,
if\footnote{Hereafter, $p$ denotes both parity defining a
superstructure and the characteristic of the ground field; the
context is, however, always clear.}
\[
 B(F x, y)+(-1)^{p(x)p(F)}B(x, Fy)=0\quad\text{for any~}x,y\in V.
\]
All linear maps preserving a given bilinear form constitute a Lie
sub(super)algebra $\faut_B(V)$ of $\fgl(V)$ denoted
$\faut_B(n)\subset\fgl(n)$ in matrix realization and consisting of
the supermatrices $X$ such that
\[
BX+(-1)^{p(X)}X^{st}B=0, \] where the {\it supertransposition} $st$
acts as follows (in the standard format):
\[
st\colon  \begin{pmatrix}A&B\\C&D\end{pmatrix}\tto
\begin{pmatrix}A^t&-C^t\\B^t&D^t\end{pmatrix}.
\]

A) The case of purely even space $V$ of dimension $n$ over a
field of characteristic $p\neq 2$. Every non-zero form $B$ can be
uniquely represented as the sum of a symmetric and an anti-symmetric
form and it is possible to consider automorphisms and equivalence
classes of each summand separately.

If the ground field $\Kee$ of characteristic $p>2$
satisfies\footnote{In this paper, $\Kee$ is algebraically closed;
over fields  algebraically non-closed, there are more types of
symmetric forms.} $\Kee^2=\Kee$, then there is just one equivalence
class of non-degenerate symmetric even forms, and the corresponding
Lie algebra $\faut_B(V)$ is called {\it orthogonal} and denoted
$\fo_B(n)$ (or just $\fo(n)$). Non-degenerate anti-symmetric forms
over $V$ exist only if $n$ is even; in this case, there is also just
one equivalence class of non-degenerate antisymmetric even forms;
the corresponding Lie algebra $\faut_B(n)$ is called  {\it
symplectic} and denoted $\fsp_B(2k)$ (or just $\fsp(2k)$). Both
algebras $\fo(n)$ and $\fsp(2k)$ are simple.

{\bf If $p=2$, the space of anti-symmetric bilinear forms is a
subspace of symmetric bilinear forms}. Also, instead of a unique
representation of a given form as a sum of an anti-symmetric and
symmetric form, we have a subspace of symmetric forms and the
quotient space of non-symmetric forms; it is not immediately clear
what to take for a representative of a given non-symmetric form. For
an answer and classification, see Lebedev's thesis \cite{LeD} and
\cite{Le1}. There are no new simple Lie superalgebras associated
with non-symmetric forms, so we confine ourselves to symmetric ones.


Instead of orthogonal and symplectic Lie algebras we have two
different types of orthogonal Lie algebras (see Theorem
\ref{SymForm}). Either the derived algebras of these algebras or
their quotient modulo center are simple if $n$ is large enough, so
the canonical expressions of the forms $B$ are needed as a step
towards classification of simple Lie algebras in characteristic $2$
which is an open problem, and as a step towards a version of this
problem for Lie superalgebras, even less investigated.

In \cite{Le1}, Lebedev showed that, with respect to the above
natural equivalence of forms (\ref{eqformM}), the following fact
takes place:

\parbegin{Theorem}[\cite{Le1}]\label{SymForm}\label{s2.2.1} Let $\Kee$ be a
perfect (i.e., such that every element of $\Kee$ has a square
root\footnote{Since $a^2-b^2=(a-b)^2$ if $p=2$, it follows that no
element can have two distinct square roots.}) field of
characteristic $2$. Let $V$ be an $n$-di\-men\-sio\-nal space over
$\Kee$.

\textup{1)} For $n$ odd, there is only one equivalence class of
non-degenerate symmetric bilinear forms on $V$.

\textup{2)} For $n$ even, there are two equivalence classes of
non-degenerate symmetric bilinear forms, one --- with at least one
non-zero element on the main diagonal of its Gram matrix ---
contains $1_n$ and the other one --- all its Gram matrices are {\it
zero-diagonal}
--- contains $S_n:=\antidiag(1, \dots, 1)$ and $\Pi_n$, where
\[\Pi_n=\begin{cases} \mat{0&1_k\\
1_k&0}&\text{if $n=2k$},\\[3mm]
\mat{0&0&1_k\\
0&1&0\\
1_k&0&0}&\text{if $n=2k+1$}.\end{cases}
\]
\end{Theorem}

Thus, every {\bf even} symmetric non-degenerate form on a superspace
of dimension $n_{\ev}|n_{\od}$ over $\Kee$ is equivalent to a form
of the shape (here: $i=\bar 0$ or $\bar 1$ and each $n_i$ may equal
to 0),
\[
B=\mat{ B_{\ev}&0\\0&B_{\od}}, \quad \text{where
$B_i=\begin{cases}1_{n_i}&\text{if $n_i$ is odd,}\\
\text{either $1_{n_i}$ or $\Pi_{n_i}$}&\text{if $n_i$ is
even.}\end{cases}$}
\]
In other words, the bilinear forms with matrices $1_{n}$ and
$\Pi_{n}$ are equivalent if $n$ is odd and non-equivalent if $n$ is
even. The Lie superalgebra preserving the bilinear form $B$
is spanned by the supermatrices which in the standard format are
of the form
\[
\mat{ A_{\ev}&B_{\ev}C^TB_{\od}^{-1}\\C&A_{\od}}, \
\begin{matrix}\text{where
$A_{\ev}\in\fo_{B_{\ev}}(n_\ev)$, $A_{\od}\in\fo_{B_{\od}}(n_\od)$, and}\\
\text{$C$ is arbitrary $n_{\od}\times n_{\ev}$ matrix.}\end{matrix}
\]
By analogy with the orthosymplectic Lie superalgebras $\fosp$ in
characteristic $0$ we call the Lie superalgebra preserving the bilinear form $B$
{\it ortho-orthogonal} and denote
$\foo_B(n_\ev|n_\od)$; usually, for clarity, we denote it
$\foo_{B_{\ev}B_{\od}}(n_\ev|n_\od)$, in particular, if $B_{\ev}=1_{n_\ev}$
and $B_{\od}=\Pi_{n_\od}$ we write $\foo_{I\Pi}(n_\ev|n_\od)$.

Since, as is easy to see, \[\foo_{\Pi I}(n_\ev|n_\od)\simeq
\foo_{I\Pi}(n_\od|n_\ev),\] we do not have to consider the Lie
superalgebra $\foo_{\Pi I}(n_\ev|n_\od)$ separately in many questions,
unless we study
Cartan prolongations where the difference between these two
incarnations of one algebra is vital: For the one, the prolong is
finite dimensional (the automorphism algebra of the $p=2$ analog of
the Riemann geometry), for the other one it is infinite dimensional
(an analog of the Lie superalgebra of Hamiltonian vector fields).

B) For an {\bf odd} symmetric form $B$ on a superspace of dimension
$(n_{\ev}|n_{\od})$ to be
non-degenerate, we need $n_{\ev}=n_{\od}$, and every such form $B$
is equivalent to $\Pi_{k|k}$, where $k=n_{\ev}=n_{\od}$, and which
is same as $\Pi_{2k}$ if the superstructure is forgotten. This form
is preserved, over $\Kee$ for $\Char\Kee\neq 2$, by linear transformations with supermatrices in the
standard format of the shape
\begin{equation}
\label{pe} \mat{ A&C\\D&A^T}, \quad \text{where $A\in\fgl(k)$, $C=C^T$
and $D=-D^T$}. 
\end{equation}
The Lie superalgebra of linear maps
preserving $B$ will be referred to as {\it periplectic}, as A.~Weil
suggested, and denoted $\fpe_B(k)$ or just $\fpe(k)$.

Note that even
the superdimensions of the characteristic $2$ versions of the Lie
(super)algebras $\faut_B(k)$ differ from their analogs in other
characteristics for both even and odd forms $B$.

C) Observe that 
\begin{equation}\label{nt}
\begin{tabular}{c}
{\bf The fact that two bilinear forms are inequivalent does}\\
{\bf not, generally, imply that the Lie (super)algebras that}\\
{\bf preserve them are not isomorphic}.
\end{tabular}
\end{equation} In \cite{Le1}, Lebedev proved that for the {\it non-degenerate
symmetric} forms, the implication spoken about in (\ref{nt}) is,
however, true (bar a few exceptions), and therefore we have several
types of non-isomorphic Lie (super) algebras (except for occasional
isomorphisms intermixing the types, e.g.,
$\foo_{I\Pi}\simeq\foo_{\Pi I}$ and
$\foo_{\Pi\Pi}'(6|2)\simeq\fpe'(4)$).

The problem of describing preserved bilinear forms has two levels:
we can consider {\it linear transformations} (Linear Algebra) and
{\it arbitrary coordinate changes} (Differential Geometry). 
In the literature, both levels are completely investigated, except
for the case where $p=2$. More precisely, the fact that the
non-split and split forms of the Lie algebras that preserve the
symmetric bilinear forms are not always isomorphic was never
mentioned. (Although known for the Chevalley groups preserving these
forms, cf. \cite{St}, these facts do not follow from each other
since there is no analog of Lie theorem on the correspondence
between Lie groups and Lie algebras.) Here we consider the Linear
Algebra aspect, for the Differential Geometry related to the objects
considered here, see \cite{LeP}.


\paragraph{Known facts: The case $p= 2$}\label{s2.1.5} The following facts
are given for clarity: lecturing on these results during the past
several years we have encountered incredulity of the listeners based
on several false premises intermixed with correct statements.

With any symmetric bilinear form $B$ the quadratic form $Q(x):=B(x,
x)$ is associated. Arf has discovered {\it the Arf invariant} --- an
important invariant of non-degenerate quadratic forms in
characteristic $2$. Two such forms are equivalent if and only if
their Arf invariants are equal, see \cite{Dye}.

The other way round, given a quadratic form $Q$, one defines a
symmetric bilinear form, called {\it the polar form}\index{bilinear
form! polar} of $Q$, by setting
\[
B_Q(x, y)=Q(x+y)-Q(x)-Q(y).
\]
The Arf invariant can not, however, be used for classification of
symmetric bilinear forms because one symmetric bilinear form can
serve as the polar form for two non-equivalent (and having different
Arf invariants) quadratic forms. Moreover, {\bf not every symmetric
bilinear form can be represented as a polar form. If $p= 2$, the
correspondence $Q\longleftrightarrow B_Q$ is not one-to-one}.

In view of (\ref{nt}) the statement of the next Lemma (proved in
\cite{Le1}) is non-trivial.

\parbegin{Lemma}\label{noniso} \textup{1)} The Lie algebras
$\fo_I(2k)$ and $\fo_\Pi(2k)$ are not isomorphic (though are of the
same dimension); the same applies to their derived algebras:

\textup{2)} $\fo_I'(2k)\not\simeq \fo_\Pi'(2k)$, though
$\dim\fo_I'(2k)=\dim\fo_\Pi'(2k)$;

\textup{3)} $\fo_I^{(2)}(2k)\not\simeq \fo_\Pi^{(2)}(2k)$ unless
$k=1$.
\end{Lemma}

Based on these results, Lebedev described all the (five) possible
analogs of the Poisson bracket, and (there exists just one) contact
bracket. Similar results for the odd bilinear form yield a
description of the anti-bracket (a.k.a. Schouten or Buttin bracket),
and the (peri)contact bracket, compare \cite{LeP} with \cite{LSh}.
The quotients of the Poisson and Buttin Lie (super)algeb\-ras modulo
center --- analogs of Lie algebras of Hamiltonian vector fields, and
their divergence-free subalgebras
--- are also described in \cite{LeP}.

\ssec{Analogs of functions and vector fields for $p>0$}
\sssec{Divided powers} Let us consider the supercommutative
superalgebra $\Cee[x]$ of polynomials in $a$ indeterminates $x =
(x_1,...,x_a)$, for convenience ordered in a ``standard format'',
i.e., so that the first $m$ indeterminates are even and the rest $n$
ones are odd ($m+n=a$). Among the integer bases of $\Cee[x]$ (i.e.,
the bases, in which the structure constants are integers), there are
two canonical ones,
--- the usual, monomial, one and the basis of {\it divided powers},
\index{divided power} which is constructed in the following way.

For any multi-index $\underline{r}=(r_1, \ldots , r_a)$, where
$r_1,\dots,r_m$ are non-negative integers, and $r_{m+1},\dots,r_n$
are $0$ or $1$, we set
\[
u_i^{(r_{i})} := \frac{x_i^{r_{i}}}{r_i!}\quad \text{and}\quad
u^{(\underline{r})} := \prod\limits_{i=1}^a u_i^{(r_{i})}.
\]
These $u^{(\underline{r})}$ form an integer basis of $\Cee[x]$.
Clearly, their multiplication relations are
\begin{equation}
\label{divp}
\renewcommand{\arraystretch}{1.4}
\begin{array}{l}  u^{(\underline{r})} \cdot u^{(\underline{s})} =
\prod\limits_{i=m+1}^n
\min(1,2-r_i-s_i)\cdot(-1)^{\sum\limits_{m<i<j\leq a} r_js_i}\cdot
\binom {\underline{r} + \underline{s}} {\underline{r}}
u^{(\underline{r} + \underline{s})},  \\
\text{where}\quad \binom {\underline{r} + \underline{s}}
{\underline{r}}:=\prod\limits_{i=1}^m\binom {r_{i} + s_{i}} {r_{i}}.
\end{array}
\end{equation}
In what follows, for clarity, we will write exponents of divided
powers in parentheses, as above, especially if the usual exponents
might be encountered as well.

Now, for an arbitrary field $\Kee$ of characteristic $p>0$, we may
consider the supercommutative superalgebra $\Kee[u]$ spanned by
elements $u^{(\underline{r})}$ with multiplication relations
(\ref{divp}). For any $m$-tuple $\un = (N_1,..., N_m)$, where $N_i$
are either positive integers or infinity, denote (we assume that
$p^\infty=\infty$)
\begin{equation}
\label{u;N} \cO(m; \un):=\Kee[u;
\un]:=\Span_{\Kee}\left(u^{(\underline{r})}\mid r_i
\begin{cases}< p^{N_{i}}&\text{for $i\leq m$}\\
=0\text{ or 1}&\text{for $i>m$}\end{cases}\right).
\end{equation}
From \eqref{divp} it is clear that $\Kee[u; \un]$ is a subalgebra of
$\Kee[u]$.  The algebra $\Kee[u]$ and its subalgebras $\Kee[u; \un]$
are called the {\it algebras of divided powers;} they can be
considered as analogs of the polynomial algebra. An important
particular case:  $\cO(m; \un_s):=\Kee[u; \un_s]$, where
$\un_s:=(1,\dots,1)$, is the algebra of truncated polynomials.

Only one of these numerous algebras of divided powers $\cO(n;\un)$
are indeed generated by the indeterminates declared: If $N_i=1$ for
all $i$. Otherwise, in addition to the $u_i$, we have to add
$u_i^{(p^{k_i})}$ for all $i\leq m$ and all $k_i$ such that
$1<k_i<N_i$ to the list of generators. Since any derivation $D$ of a
given algebra is determined by the values of $D$ on the generators,
we see that $\fder(\cO[m; \un])$ has more than $m$ functional
parameters (coefficients of the analogs of partial derivatives) if
$N_i\neq 1$ for at least one $i$. Define {\it distinguished partial
derivatives}\index{derivative! partial, distinguished}
 by setting
\[
\partial_i(u_j^{(k)})=\delta_{ij}u_j^{(k-1)}\ \text{ for any $k<p^{N_j}$}.
\]

The simple vectorial Lie algebras over $\Cee$ have only one
parameter: the number of indeterminates. If $\Char~ \Kee =p>0$, the
vectorial Lie algebras acquire one more parameter: $\un$. For Lie
superalgebras, $\un$ only concerns the even indeterminates.

The Lie (super)algebra of all derivations $\fder(\cO[m; \un])$ turns
out to be not so interesting as its {\it Lie subsuperalgebra of
distinguished derivations}: Let
\begin{equation}\label{vect-super}
\begin{array}{c} \fvect (m; \un|n) \ \text{ a.k.a. }W(m;
\un|n)\  \text{ a.k.a. }\\
\fder_{dist} \Kee[u;
\un]=\Span_{\Kee}\left(u^{(\underline{r})}\partial_k\mid
r_i\begin{cases}< p^{N_{i}}&\text{for $i\leq m$},\\
=0\text{ or 1}&\text{for $i>m$};\end{cases} \quad  1\leq k\leq
n\right)\end{array}
\end{equation} be the {\it general vectorial Lie algebra of
distinguished derivations}. The next notions are analogs of the
polynomial algebra of the dual space.

\sssec{Recapitulation: On vectorial Lie superalgebras, there are TWO
analogs of trace}\label{diverg} More precisely, there are {\it
traces} and their Cartan prolongs, called {\it divergencies}. On any
Lie (super)algebra $\fg$ over a field $\Kee$, a {\it trace} is any
map $\tr: \fg\tto \Kee$ such that
\begin{equation}
\label{deftr} \tr ([\fg,\fg])=0.
\end{equation}

The straightforward analogs of the trace are, therefore, the linear
functionals that vanish on $\fg':=[\fg,\fg]$; the number of linearly
independent traces is equal to $\dim \fg/\fg'$; if $\fg$ is a Lie
superalgebra, these traces are called supertraces and they can be
even or odd. Each trace is defined up to a non-zero scalar factor
selected {\it ad lib}.

Let now $\fg$ be a $\Zee$-graded vectorial Lie superalgebra with
$\fg_{-}:=\mathop{\oplus}\limits_{i<0}\fg_i$ generated by
$\fg_{-1}$, and let $\tr$ be a (super)trace on $\fg_0$. The {\it
divergence} $\Div:\fg\tto\cF$, where $\cF$ is the space of
functions-coefficients, is an $\ad_{\fg_{-1}}$-invariant
prolongation of the trace satisfying the following conditions:
\[
\begin{array}{l}
\Div:\fg\tto\cF\text{~~ preserves the degree, i.e., $\deg\Div=0$};\\
X_i(\Div D)=\Div[X_i,D]\text{~~ for all elements
$X_i$ that span $\fg_{-1}$};\\
 \Div|_{\fg_0}=\tr;\\
 \Div|_{\fg_{-}}=0.\end{array}\]
By construction, the Lie (super)algebra $\fs\fg:=\Ker\Div|_{\fg}$ of
divergence-free elements of $\fg$ is the complete prolong of
$(\fg_{-}, \Ker\tr|_{\fg_0})$. This fact explains why we say that
$\Div$ is the prolongation of the trace.

Strictly speaking, divergences are not traces (they do not satisfy
\eqref{deftr}) but for vectorial Lie (super)algebras they embody the
idea of the trace (understood as property \eqref{deftr}) better than
the traces. We denote the {\it special} ({\it divergence free})
subalgebra of a~vectorial algebra $\fg$ by $\fs\fg$, e.g.,
$\fsvect(n|m)$. If there are several traces on $\fg_0$, there are
several types of special subalgebras of $\fg$ and we need a
different name for each.

\ssec{Weisfeiler filtrations and gradings} Recall, see \cite{LSh},
that the \emph{Weisfeiler filtrations} were initially used for
description of simple (or primitive) transitive infinite dimensional
Lie (super)algebras $\cL$ by selecting a maximal subalgebra
$\cL_{0}$ of finite codimension. For the same reason we need these
filtrations and associated gradings dealing with infinite
dimensional algebras (if $\un_i=\infty$ for at least one $i$).

Dealing with finite dimensional algebras, we can confine ourselves
to maximal subalgebras of \emph{least} codimension, or almost least,
etc. Let $\cL_{-1}$ be a~minimal $\cL_{0}$\defis invariant subspace
strictly containing $\cL_{0}$, and $\cL_{0}$\defis invariant; for
$i\geq 1$, set:
\begin{equation}
\label{1.3} \cL_{-i-1}=[\cL_{-1}, \cL_{-i}]+\cL_{-i}\  \text{ and }\
\cL_i =\{D\in \cL_{i-1}\mid [D,
\cL_{-1}]\subset\cL_{i-1}\}.
\end{equation}
We thus get a filtration:
\begin{equation}
\label{1.2} \cL= \cL_{-d}\supset \cL_{-d+1}\supset \dots \supset
\cL_{0}\supset \cL_{1}\supset \dots
\end{equation}
The $d$ in (\ref{1.2}) is called the \emph{depth} of $\cL$ and of
the associated graded (the \emph{Weisfeiler graded}) Lie
superalgebra $\fg=\mathop{\oplus}\limits_{-d\leq i}\fg_i$, where
$\fg_{i}=\cL_{i}/\cL_{i+1}$.

\ssec{What $\fg(A)$ is}\label{Sg(A)}

\sssec{Warning: $\fpsl$ has no Cartan matrix. The relatives of
$\fsl$ and $\fpsl$ that have Cartan matrices}\label{warn} For the
most reasonable definition of Lie algebra with Cartan matrix over
$\Cee$, see \cite{K}. The same definition applies, practically
literally, to Lie superalgebras and to modular Lie algebras and to
modular Lie superalgebras. However, the usual sloppy practice is to
attribute Cartan matrices to (usually simple) Lie (super)algebras
none of which, strictly speaking, has a Cartan matrix!

Although it may look strange for those with non-super experience
over $\Cee$, neither the simple modular Lie algebra $\fpsl(pk)$, nor
the simple modular Lie superalgebra $\fpsl(a|pk+a)$, nor --- in
characteristic $0$ --- the simple Lie superalgebra $\fpsl(a|a)$
possesses a Cartan matrix. Their central extensions\footnote{If
$p=2$, the simple Lie (super)algebra may have more central
extensions; these centrally extended algebras are even further, so
to say, from possessing Cartan matrix.}
--- $\fsl(pk)$, the modular Lie superalgebra $\fsl(a|pk+a)$ for
characteristic $p>0$, and the Lie superalgebra $\fsl(a|a)$ for
characteristic $0$--- do not have Cartan matrix, either.

Their relatives possessing a Cartan matrix are, respectively,
$\fgl(pk)$, $\fgl(a|pk+a)$, and $\fgl(a|a)$, and for the grading
operator we take the matrix unit $E_{1,1}$.

Since all the Lie (super)algebras involved (the simple one, its
central extension, the derivation algebras thereof) are often needed
simultaneously (and only representatives of one of these types of
Lie (super)algebras are of the form $\fg(A)$), it is important to
have (preferably short and easy to remember) notation for each of
them. For example, in addition to $\fpsl$, $\fsl$, $\fp\fgl$ and
$\fgl$, we have:

{\bf for $p=2$}: $\fe(7)$ is of dimension 134, then
$\dim\fe(7)'=133$, whereas the ``simple core'' is $\fe(7)'/\fc$ of
dimension 132;

for an analog of $\fg(2)$, having no Cartan matrix, see
\cite{BGLLS}, where Shen's and Brown's descriptions are sharpened;

{\bf the orthogonal Lie algebras and their super analogs} are
considered in detail later.

In our main examples, $\sdim \fg(A)^{(i)}/\fc=d|\delta$  for a
simple Lie (super)algebra $\fg(A)^{(i)}/\fc$ whereas the notation
$D/d|\delta$ means that $\sdim \fg(A)=D|\delta$. The general formula
is
\begin{equation}\label{dims}
d=D-2(\size(A)-\rk(A))\text{~~and~~}i=\size(A)-\rk(A).
\end{equation}

\sssec{What Cartan matrix is}\label{sssGEN} Let $A=(A_{ij})$ be an $n\times n$-matrix with
elements in $\Kee$ with $\rk A=n-l$. Complete $A$ to an $(n+l)\times
n$-matrix $\begin{pmatrix} A\\ B \end{pmatrix}$ of rank $n$. (Thus,
$B$ is an $l\times n$-matrix.)

Let the elements $e_i^\pm, h_i$, where $i=1,\dots,n$, and $d_k$,
where $k=1,\dots,l$, generate a Lie superalgebra denoted $\widetilde
\fg(A, I)$,\index{$\widetilde\fg(A, I)$} where $I=(p_1, \dots
p_n)\in(\Zee/2)^n$ is a collection of parities ($p(e_i^\pm)=p_i$,
the parities of the $d_k$'s being $\ev$), free except for the
relations
\begin{equation}\label{gArel_0}
\renewcommand{\arraystretch}{1.4}
\begin{array}{l}
{}[e_{i}^+, e_{j}^-] = \delta_{ij}h_i; \quad [h_i, e_{j}^\pm]=\pm
A_{ij}e_{j}^\pm;\quad [d_k, e_{j}^\pm]=\pm
B_{kj}e_{j}^\pm;\\
{}[h_i, h_j]=[h_i, d_k]=[d_k, d_m]=0\text{~~for any $i, j, k,
m$}.\end{array}
\end{equation}
The Lie superalgebra  $\widetilde \fg(A, I)$ is $\Zee^n$-graded with
\begin{equation}\label{deg}
\renewcommand{\arraystretch}{1.4}
\begin{array}{l}
\deg e_{i}^\pm =(0, \ldots, 0, \pm 1, 0, \ldots, 0)\\
\deg h_i=\deg d_k=(0, \ldots, 0)\text{~~for any $i, k$}.\end{array}
\end{equation}

Let $\fh$ denote the linear span of the $h_i$'s and $d_k$'s. Let
$\widetilde \fg(A, I)^\pm$ denote the Lie subsuperalgebras in
$\widetilde \fg(A, I)$ generated by  $e_{1}^\pm, \ldots, e_{n}^\pm$.
Then
\[\widetilde \fg(A, I)=\widetilde \fg(A, I)^-\oplus \fh\oplus\widetilde \fg(A, I)^+,\]
where the homogeneous component of degree $(0, \ldots, 0)$ is just
$\fh$.

The Lie subsuperalgebras $\widetilde \fg(A, I)^\pm$ are homogeneous
in this $\Zee^n$-grading, and there is a
\begin{equation}\label{r}
\text{maximal homogeneous (in this $\Zee^n$-grading) ideal $\fr$
such that $\fr\cap\fh=0$.}
\end{equation} The ideal $\fr$ is just the sum of homogeneous ideals
whose homogeneous components of degree $(0, \ldots, 0)$ is trivial.

As $\rk A =n-l$, there exists an $l\times n$-matrix $T=(T_{ij})$ of
rank $l$ such that
\begin{equation}\label{TA=0}
TA=0.
\end{equation}
Let
\begin{equation}\label{central}
c_i=\mathop{\sum}\limits_{1\leq j\leq n} T_{ij}h_j, \text{~~where~~} i=1,\dots,l.
\end{equation}
Then, from the properties of the matrix $T$, we deduce that
\begin{equation}\label{central1}
\begin{tabular}{l}
a) the elements $c_i$ are linearly independent; let $\fc$ be the space they span;\\
b) the elements $c_i$ are central, because\\
$[c_i,e_j^\pm]=\pm\left(\sum\limits_{1\leq k\leq n}
T_{ik}A_{kj}\right) e_j^\pm=\pm (TA)_{ij} e_j^\pm
\stackrel{\eqref{TA=0}}{=}0$.
\end{tabular}
\end{equation}

The Lie (super)algebra $\fg(A, I)$\index{$\fg(A, I)$} is defined as
the quotient $\widetilde \fg(A, I)/\fr$ and is called the {\it Lie
(super)algebra with Cartan matrix $A$ (and parities $I$)}. Note that
this coincides with the definition in \cite{CE} of the {\it
contragredient}\footnote{This word does not seem to mean anything in
this context, and therefore this term, though often used, is ill
chosen.} Lie superalgebras, although written in a slightly different
way. Condition \eqref{r} modified as
\begin{equation}\label{rc}
\text{maximal homogeneous (in this $\Zee^n$-grading) ideal $\fs$
such that $\fs\cap\fh=\fc$}
\end{equation}
leads to what in \cite{CE} is called the {\it centerless
contragredient} Lie superalgebra, cf. \cite{Bi}.

By abuse of notation we denote by $e_i^\pm, h_i, d_k$ and $\fc$
their images in $\fg(A,I)$ and $\fg(A,I)'$.

The Lie superalgebra $\fg(A,I)$ inherits, clearly, the
$\Zee^n$-grading of $\widetilde \fg(A, I)$. The non-zero elements
$\alpha\in\Zee^n\subset \Ree^n$ such that the homogeneous component
$\fg(A,I)_\alpha$ is non-zero are called {\it roots}.
\index{Root}The set $R$ of all roots is called {\it the root
system}\index{System!root} of $\fg$. Clearly, the subspaces
$\fg_\alpha$ are purely even or purely odd, and the corresponding
roots are said to be \textit{even} or \textit{odd}.

The additional to \eqref{gArel_0} relations that turn $\widetilde
\fg(A, I)^\pm$ into $\fg(A, I)^\pm$ are of the form $R_i=0$ whose
left sides are implicitly described as follows:
\begin{equation}
\label{myst1}
\begin{split}
 &\text{the $R_i$ that generate the maximal ideal $\fr$.}
 \end{split}
\end{equation}
{\bf For the explicit description of these additional relations},
see \cite{BGLL1}.


\sssec{Roots and weights}\label{roots} In this subsection, $\fg$
denotes one of the algebras $\fg(A,I)$ or $\widetilde{\fg}(A,I)$.

The elements of $\fh^*$ are called {\it weights}.\index{weight} For
a given weight $\alpha$, the {\it weight subspace} of a given
$\fg$-module $V$ is defined as
\[
V_\alpha=\{x\in V\mid \text{an integer $N>0$ exists such that
$(\alpha(h)-\ad_h)^N x=0$ for any $h\in\fh$}\}.
\]

Any non-zero element $x\in V$ is said to be {\it of weight
$\alpha$}. For the roots, which are particular cases of weights if
$p=0$, the above definition is inconvenient: In the modular analog
of the following useful statement summation should be over roots
defined in the previous subsection.

\parbegin{Statement}[\cite{K}] Over $\Cee$, the space $\fg$ can be
represented as a direct sum of subspaces
\[
\fg=\mathop{\bigoplus}\limits_{\alpha\in \fh^*} \fg_\alpha.
\]
\end{Statement}

Note that $\fh\subsetneq\fg_0$ over $\Kee$, e.g., all weights of the
form $p\alpha$ over $\Cee$ become 0.

\sssec{Systems of simple and positive roots} In this subsection,
$\fg=\fg(A,I)$, and $R$ is the root system of $\fg$.

For any subset $B=\{\sigma_{1}, \dots, \sigma_{m}\} \subset R$, we
set  (we denote by $\Zee_{+}$ the set of non-negative integers):
\[
R_{B}^{\pm} =\{ \alpha \in R \mid \alpha = \pm \sum n_{i}
\sigma_{i},\ \ n_{i} \in \Zee_{+} \}.
\]

The set $B$ is called a {\it system of simple roots~}\index{Root!
simple system of}\  of $R$ (or $\fg$) if $ \sigma_{1}, \dots ,
\sigma_{m}$ are linearly independent and $R=R_B^+\cup R_B^-$. Note
that $R$ contains basis coordinate vectors, and therefore spans
$\Ree^n$; thus, any system of simple roots contains exactly $n$
elements.

Let $(\cdot,\cdot)$ be the standard Euclidean inner product in
$\Ree^n$. A subset $R^+\subset R$ is called a {\it system of
positive roots~}\index{Root! positive system of}\  of $R$ (or $\fg$)
if there exists $x\in\Ree^n$ such that
\begin{equation}\label{x}
\begin{split}
 &(\alpha,x)\in\Ree\backslash \{0\}\text{ for any $\alpha\in R$},\\
 &R^+=\{\alpha\in R\mid (\alpha,x)>0\}.
\end{split}
\end{equation} Since $R$ is a finite (or, at least, countable
if $\dim \fg(A)=\infty$) set, so the set
\[\{y\in\Ree^n\mid\text{there exists $\alpha\in R$ such that }
(\alpha,y)=0\}
\]
is a finite/countable union of $(n-1)$-dimensional subspaces in
$\Ree^n$, so it has zero measure. So for almost every $x$, condition
(\ref{x}) holds.

By construction, any system $B$ of simple roots is contained in
exactly one system of positive roots, which is precisely $R_B^+$.

\parbegin{Statement} Any finite system $R^+$ of positive roots of
$\fg$ contains exactly one system of simple roots. This system
consists of all the positive roots (i.e., elements of $R^+$) that
can not be represented as a sum of two positive
roots.\end{Statement}

We can not give an {\it a priori} proof of the fact that each set of
all positive roots each of which is not a sum of two other positive
roots consists of linearly independent elements. This is, however,
true for finite dimensional Lie algebras and superalgebras $\fg(A,
I)$ if $p\neq 2$.

\sssec{Normalization convention}\label{normA} Clearly,
\begin{equation}
\label{rescale} \text{the rescaling
$e_i^\pm\mapsto\sqrt{\lambda_i}e_i^\pm$, sends $A$ to $A':=
\diag(\lambda_1, \dots , \lambda_n)\cdot A$.} 
\end{equation}
Two pairs $(A, I)$ and $(A', I')$ are said to be {\it equivalent}
(and we write $(A, I)\sim(A', I')$) if $(A', I')$ is obtained from
$(A, I)$ by a composition of a permutation of parities and a
rescaling $A' = \diag (\lambda_{1}, \dots, \lambda_{n})\cdot A$,
where $\lambda_{1}\dots \lambda_{n}\neq 0$. Clearly, equivalent
pairs determine isomorphic Lie superalgebras.

The rescaling affects only the matrix $A_B$, not the set of parities
$I_B$. The Cartan matrix $A$ is said to be {\it
normalized}\index{Cartan matrix, normalized} if
\begin{equation}
\label{norm} A_{jj}=0\text{~~ or 1, or 2,}
\end{equation}
where we let $A_{jj}=2$ only if $p_j=\ev$; in order to distinguish
between the cases where $p_j=\ev$ and $p_j=\od$, we write
$A_{jj}=\ev$ or $\od$, instead of 0 or 1, if $p_j=\ev$. {\bf We will
only consider normalized Cartan matrices; for them, we do not have
to describe $I$.}

The row with a $0$ or $\ev$ on the main diagonal can be multiplied
by any nonzero factor; usually (not only in this paper) we multiply
the rows so as to make $A_{B}$ symmetric, if possible.

{\it A posteriori}, for each {\bf finite dimensional} Lie
(super)algebra of the form $\fg(A)$ with indecomposable Cartan
matrix $A$, the matrix $A$ is symmetrizable (i.e., it can be made
symmetric by operation (\ref{rescale})) for any $p$. For affine and
almost affine Lie (super)algebra of the form $\fg(A)$ this is not
so, cf. \cite{CCLL}.

\sssec{Equivalent systems of simple roots} \label{EqSSR} Let
$B=\{\alpha_1,\dots,\alpha_n\}$ be a system of simple roots. Choose
non-zero elements $e_i^\pm$ in the 1-dimensional (by definition)
superspaces $\fg_{\pm\alpha_i}$; set $h_{i}=[e_{i}^{+}, e_{i}^-]$,
let $A_{B} =(A_{ij})$, where the entries $A_{ij}$ are recovered from
relations \eqref{gArel_0}, and let $I_{B}=\{p(e_{1}), \cdots,
p(e_{n})\}$. Lemma \ref{serg} claims that all the pairs $(A_B,I_B)$
are equivalent to each other.

Two systems of simple roots $B_{1}$ and $B_{2}$ are said to be {\it
equivalent} if the pairs $(A_{B_{1}}, I_{B_{1}})\sim(A_{B_{2}},
I_{B_{2}})$.

For the role of the ``best'' (first among equals) order of indices
we propose the one that minimizes the value
\begin{equation}\label{minCM}
\max\limits_{i,j\in\{1,\dots,n\}\text{~such that~}(A_B)_{ij}\neq
0}|i-j|
\end{equation}
(i.e., gather the non-zero entries of $A$ as close to the main
diagonal as possible).

\sssec{Chevalley generators and Chevalley bases}\label{SsChev} We
often denote the set of generators corresponding to a normalized
matrix by $X_{1}^{\pm},\dots , X_{n}^{\pm}$ instead of
$e_{1}^{\pm},\dots , e_{n}^{\pm}$; and call them, together with the
elements $H_i:=[X_{i}^{+}, X_{i}^{-}]$, and the derivatives $d_j$
added for convenience for all $i$ and $j$, the {\it Chevalley
generators}.\index{Chevalley generator}

For $p=0$ and normalized Cartan matrices of simple finite
dimensional Lie algebras, there exists only one (up to signs) basis
containing $X_i^\pm$ and $H_i$  in which $A_{ii}=2$ for all $i$ and
all structure constants are integer, cf. \cite{St}. Such a basis is
called the {\it Chevalley}\index{Basis! Chevalley}  basis.

Observe that, having normalized the Cartan matrix of $\fo(2n+1)$ so
that $A_{ii}=2$ for all $i\neq n$ but $A_{nn}=1$, we get {\bf
another} basis with integer structure constants. Clearly, this basis
also qualifies to be called {\it Chevalley  basis}; for the Lie
superalgebras, the basis normalized as in \eqref{norm} is more
appropriate than the one with $A_{ii}=2$; for $p=2$, the
normalization \eqref{norm} seems at the moment the only reasonable
one:

\sssbegin{Conjecture} If $p>2$, then for finite dimensional Lie
(super)algebras with indecomposable Cartan matrices normalized as in
$(\ref{norm})$, there also exists only one (up to signs) analog of
the Chevalley basis.

The analogs of Chevalley bases for $p=2$ are not described yet; we
conjecture that the methods of a recent paper \cite{CR} should solve
the problem.
\end{Conjecture}

\ssec{Ortho-orthogonal and periplectic Lie superalgebras}\label{Soo}

In this section, $p=2$ and $\Kee$ is perfect. We also assume that
$n_\ev,n_\od>0$.

\sssec{Non-degenerate even supersymmetric bilinear forms and
ortho\Defis or\-tho\-gonal Lie superalgebras} For $p=2$, there are,
in general, four equivalence classes of inequivalent non-degenerate
even supersymmetric bilinear forms on a given superspace. Any such
form $B$ on a superspace $V$ of superdimension $n_\ev|n_\od$ can be
decomposed as follows:
\[
B=B_\ev\oplus B_\od,
\]
where $B_\ev$, $B_\od$ are symmetric non-degenerate forms on $V_\ev$
and $V_\od$, respectively. For $i=\ev,\od$, the form $B_i$ is
equivalent to $1_{n_i}$ if $n_i$ is odd, and equivalent to $1_{n_i}$
or $\Pi_{n_i}$ if $n_i$ is even. So every non-degenerate even
symmetric bilinear form is equivalent to one of the following forms
(some of them are defined not for all dimensions):
\[
\begin{array}{ll}
B_{II}=1_{n_\ev}\oplus 1_{n_\od}; &B_{I\Pi}=1_{n_\ev}\oplus
\Pi_{n_\od}\text{~if~}n_\od \text{~is even;}\\
B_{\Pi I}=\Pi_{n_\ev}\oplus 1_{n_\od}\text{~if~}n_\ev \text{~is
even}; &B_{\Pi\Pi}=\Pi_{n_\ev}\oplus \Pi_{n_\od}\text{~if~}n_\ev,
n_\od \text{~are even.}
\end{array}
\]
We denote the Lie superalgebras that preserve the respective forms
by $\fo\fo_{II}(n_\ev|n_\od)$, $\fo\fo_{I\Pi}(n_\ev|n_\od)$,
$\fo\fo_{\Pi I}(n_\ev|n_\od)$, $\fo\fo_{\Pi\Pi}(n_\ev|n_\od)$,
respectively. Now let us describe these algebras.

\paragraph{$\fo\fo_{II}(n_\ev|n_\od)$} If $n\geq 3$, then the Lie
superalgebra $\fo\fo_{II}'(n_\ev|n_\od)$ is simple. This Lie
superalgebra {\bf has no Cartan matrix}.

\paragraph{$\fo\fo_{I\Pi}(n_\ev|n_\od)$~~($n_\od=2k_\od$)} The Lie
superalgebra $\fo\fo_{I\Pi}'(n_\ev|n_\od)$ is simple, it has Cartan
matrix if and only if $n_\ev$ is odd; this matrix has the following
form (up to a format; all possible formats --- corresponding to
$\ast=0$ or $\ast=\ev$ --- are described in Table \ref{tbl} below):
\begin{equation}\label{oowith1}
\begin{pmatrix} \ddots&\ddots&\ddots&\vdots\\
\ddots&\ast&1&0\\
\ddots&1&\ast&1\\
\cdots&0&1&1\end{pmatrix}
\end{equation}
In particular, the Lie algebra $\fg=\fo^{(1)}(2n+1)$ with Cartan
matrix \eqref{oowith1} with $\ast=\ev$ can be considered as the Lie
algebra of matrices of the form (recall that $ZD(n)$ is the space of
symmetric matrices with zeros on the main diagonal)
\begin{equation}\label{oowith1matr}
\begin{pmatrix}A&X&B\\Y^T&0&X^T\\C&Y&A^T\end{pmatrix},\text{~where~}\begin{array}{c}
A\in\fgl(n);~~B,C\in ZD(n);\\X,Y\text{~are column
$n$-vectors.}\end{array}
\end{equation}

\paragraph{$\fo\fo_{\Pi\Pi}(n_\ev|n_\od)~~(n_\ev=2k_\ev,
n_\od=2k_\od)$}\label{oo-oo_PP} If $n=n_\ev+n_\od\geq 6$, then
\begin{equation}
\label{oopipi}
\begin{split}
&\text{if $k_\ev+k_\od$ is odd, then the Lie superalgebra
$\fo\fo_{\Pi\Pi}^{(2)}(n_\ev|n_\od)$ is simple;} \cr &\text{if
$k_\ev+k_\od$ is even, then the Lie superalgebra
$\fo\fo_{\Pi\Pi}^{(2)}(n_\ev|n_\od)/\Kee 1_{n_\ev|n_\od}$ is
simple.}
\end{split}
\end{equation}

Each of these simple Lie superalgebras is also close to a Lie
superalgebra with Cartan matrix. To describe this Cartan matrix Lie
superalgebra in most simple terms, we will choose a slightly
different realization of $\fo\fo_{\Pi\Pi}(2k_\ev|2k_\od)$: Let us
consider it as the algebra of linear transformations that preserve
the bilinear form $\Pi_{2k_\ev+2k_\od}$ in the supermatrix format
$k_\ev|k_\od|k_\ev|k_\od$. Then the algebra
$\fo\fo_{\Pi\Pi}^{(i)}(2k_\ev|2k_\od)$ is spanned by supermatrices
of format $k_\ev|k_\od|k_\ev|k_\od$ and of the form
\begin{equation}\label{matform}
\begin{pmatrix}A&C\\D&A^T\end{pmatrix}\text{~where~}
\begin{array}{l}
A\in\begin{cases}\fgl(k_\ev|k_\od)&\text{if~}i\leq
1,\\\fsl(k_\ev|k_\od)&\text{if~}i\geq 2,\end{cases} \\ 
C,D\text{~are~}\begin{cases}\text{symmetric matrices}&\text{if~}
i=0;\\\text{symmetric zero-diagonal matrices}&\text{if~} i\geq
1.\end{cases}\end{array}
\end{equation}
If $i\geq 1$, these derived algebras have a non-trivial central
extension given by the following cocycle:
\begin{equation}\label{cocycle}
F\left(\begin{pmatrix}A&C\\D&A^T\end{pmatrix},
\begin{pmatrix}A'&C'\\D'&A'^T\end{pmatrix}\right)=\sum\limits_{1\leq
i<j\leq k_\ev+k_\od} (C_{ij}D'_{ij}+C'_{ij}D_{ij})
\end{equation}
(note that this expression resembles $\frac 12\tr(CD'+C'D)$). We
will denote this central extension of
$\fo\fo_{\Pi\Pi}^{(i)}(2k_\ev|2k_\od)$ by
$\fo\fo\fc(i,2k_\ev|2k_\od)$.

Let\index{$I_0:=\diag(1_{k_\ev\vert k_\od},0_{k_\ev\vert k_\od})$}
\begin{equation}
\label{I_0osp}I_0:=\diag(1_{k_\ev|k_\od},0_{k_\ev|k_\od}).
\end{equation}
Then the corresponding Cartan matrix Lie superalgebra is
\begin{equation}
\label{ooc}
\begin{split}
&\fo\fo\fc(2,2k_\ev|2k_\od)\ltimes\Kee I_0\quad\text{ if
$k_\ev+k_\od$ is odd;} \cr &\fo\fo\fc(1,2k_\ev|2k_\od)\ltimes\Kee
I_0\quad\text{ if $k_\ev+k_\od$ is even.}
\end{split}
\end{equation}

The corresponding Cartan matrix has the following form (up to a
format; all possible formats --- corresponding to $\ast=0$ or
$\ast=\ev$ --- are described in Table \ref{tbl} below):
\begin{equation}\label{ooPPCM}
\begin{pmatrix}
\ddots&\ddots&\ddots&\vdots&\vdots\\
\ddots&\ast&1&0&0\\
\ddots&1&\ast&1&1\\
\cdots&0&1&\ev&0\\
\cdots&0&1&0&\ev\end{pmatrix}
\end{equation}

\sssec{The non-degenerate odd supersymmetric bilinear forms.
Periplectic Lie superalgebras} \label{peLS} In this subsection, $m\geq
3$.

\begin{equation}
\label{pe0}
\begin{split}
&\text{If $m$ is odd, then the Lie superalgebra $\fpe_B^{(2)}(m)$ is
simple;} \cr &\text{If $m$ is even, then the Lie superalgebra
$\fpe_B^{(2)}(m)/\Kee 1_{m|m}$ is simple.}
\end{split}
\end{equation}

If we choose the form $B$ to be $\Pi_{m|m}$, then the algebras
$\fpe_B^{(i)}(m)$ consist of matrices of the form (\ref{matform});
the only difference from $\fo\fo_{\Pi\Pi}^{(i)}$ is the format which
in this case is $m|m$.

Each of these simple Lie superalgebras has a $2$-structure. Note
that if $p\neq 2$, then the Lie superalgebra $\fpe_B(m)$ and its
derived algebras are not close to Cartan matrix Lie superalgebras
(because, for example, their root system is not symmetric). If $p=2$
and $m\geq 3$, then they {\bf are} close to Cartan matrix Lie
superalgebras; here we describe them.

The algebras $\fpe_B^{(i)}(m)$, where $i>0$, have non-trivial
central extensions with cocycles (\ref{cocycle}); we denote these
central extensions by $\fpe\fc(i,m)$. Let us introduce another
matrix \index{$I_0:=\diag(1_m,0_m)$}
\begin{equation}
\label{I_0pe}I_0:=\diag(1_m,0_m).
\end{equation} Then the Cartan matrix Lie
superalgebras are
\begin{equation}
\label{pec}
\begin{split}
&\fpe\fc(2,m)\ltimes\Kee I_0\text{ if $m$ is odd;} \cr
&\fpe\fc(1,m)\ltimes\Kee I_0\text{ if $m$ is even.}
\end{split}
\end{equation}

The corresponding Cartan matrix has the form (\ref{ooPPCM}); the
only condition on its format is that the last two simple roots must
have distinct parities. The corresponding Dynkin diagram is shown in
Table \ref{tbl}; all its nodes, except for the ``horns'', may be
both $\otimes$ or~$\odot$, see (\ref{cm1}).

\sssec{Superdimensions} The following expressions (with a $+$ sign)
are the superdimensions of the
 relatives of the ortho-orthogonal and
periplectic Lie superalgebras that possess Cartan matrices. To get
the superdimensions of the simple relatives, one should replace $+2$
and $+1$ by $-2$ and $-1$, respectively, in the two first lines and
the four last ones:
\begin{equation}\label{dimtable}
\begin{array}{lll}
\dim \fo\fc (1;2k)\ltimes\Kee I_0&=2k^2-k\pm 2&\text{if $k$ is even;}\\
\dim \fo\fc (2;2k)\ltimes\Kee I_0&=2k^2-k\pm 1&\text{if $k$ is odd;}\\
\dim\fo'(2k+1)&=2k^2+k&\\
\sdim\fo\fo'(2k_\ev+1|2k_\od)&=2k_\ev^2+k_\ev+ 2k_\od^2+k_\od\mid
2k_\od(2k_\ev+1)&\\
\sdim \fo\fo\fc (1;2k_\ev|2k_\od)\ltimes\Kee I_0&=2k_\ev^2-k_\ev+
2k_\od^2-k_\od\pm 2\mid 4k_\ev k_\od&\text{if $k_\ev+k_\od$ is even;}\\
\sdim \fo\fo\fc (2;2k_\ev|2k_\od)\ltimes\Kee I_0&=2k_\ev^2-k_\ev+
2k_\od^2-k_\od\pm 1\mid 4k_\ev k_\od&\text{if $k_\ev+k_\od$ is odd;}\\
\sdim \fpe\fc (1;m)\ltimes\Kee I_0&=m^2\pm 2\mid m^2-m&\text{if $m$ is even;}\\
\sdim \fpe\fc (2;m)\ltimes\Kee I_0&=m^2\pm 1\mid m^2-m&\text{if $m$
is odd}
\end{array}
\end{equation}

\sssec{An example} Let us explain why the simple Lie algebras like
$\fpsl(np)$ over $\Kee$ of characteristic $p>0$ does not have Cartan
matrix and how its ``too small" toral subalgebra  (i.e., a
subalgebra of diagonal matrices) should be fixed (enlarged so that
the enlarged algebra would possess a Cartan matrix.

Consider the case of orthogonal Lie algebras as most complicated
one. Let $\size A=k$, i.e., consider orthogonal $2k\times
2k$-matrices. The Chevalley generators are:
\begin{equation}\label{dimtable1}
\begin{array}{ll}
e_i^+ = E^{i,i+1} + E^{k+i+1,k+i}&\text{for $i=1,...,k-1$};\\
e_k^+ = E^{k-1,2k} + E^{k,2k-1}\\
e_i^- = (e_i^+)^T&\text{for $i=1,...,k$}.
\end{array}
\end{equation}
 Let us start with $k=2n+1$ and the
algebra $\fo^{(2)}_\Pi(4n+2)$. The Cartan matrix has rank $k-1$ in
this case; the degeneration is caused by the fact that two last rows
are the same. This means that the element $h_{k-1}-h_k$ is central
in $\fg(A)$; also, this element belongs to $\fg'(A)$ since
$h_i=[e_i^+,e_i^-]$. But in the orthogonal algebra we have
\[
[e_{k-1}^+,e_{k-1}^-]=[e_k^+,e_k^-]=E^{k-1,k-1}+E^{k,k}+E^{2k-1,2k-1}+E^{2k,2k}.
\]
So we essentially have $h_{k-1}=h_k$ in the ``non-fixed" algebra. To
fix this, we need to construct a non-trivial central extension of
$\fo^{(2)}_\Pi(4n+2)$ such that $[e_{k-1}^+,e_{k-1}^-] -
[e_k^+,e_k^-]$ is the central element. The extension $\foc(2;4n+2)$
satisfies this property.

The Lie algebra $\fg(A)$ also contains an additional grading element
$d_1$ such that its action is determined by a row we add to $A$ for
it to have rank $k$. We can choose $(0,...,0,1)$ as such a row,
i.e., we have
\[
[d_1,e_i^\pm]=0\text{~~ for $i=1,...,k-1$; \quad
$[d_1,e_k^\pm]=e_k^\pm$.} \] It is easy to check that the matrix
$I_0 = \diag(1_k,0_k)$ acts in exactly this way. So $\fg(A)$ is
isomorphic to $\foc(2;4n+2) \ltimes \Kee I_0$.

Now let us consider the case $k=2n$. Let us start with
$\fo_\Pi^{(2)}(4n)$ again (not $\fo'_\Pi(4n)$). In this case the
matrix $A$ has rank $k-2$. One degeneration is again two last rows
being equal; the other one is that the sum of all odd-numbered rows
is equal to $0$. So again, first we move from $\fo_\Pi^{(2)}(4n)$ to
its central extension $\foc(2,4n)$. Fortunately, we do not need to
add another central element, the corresponding sum of $h_i$ is
already central in the algebra: $\mathop{\sum}\limits_{1\leq i\leq
n} [e_{2i-1}^+,e_{2i-1}^-]= 1_{4n}$. Now we need to add two grading
elements determined by two rows we add to $A$ to make the rank of
the enlarged matrix equal to $k$. We can choose the first row to be
$(0,...,0,1)$ again, $d_1$ is $I_0$ again. We can choose
$(1,0,...,0)$ as the second row, and the needed action coincides
with the action of the matrix $E^{1,1}+E^{k+1,k+1}$. This is one of
the matrices present in $\fo'_\Pi(4n)$ but absent in
$\fo^{(2)}_\Pi(4n)$ (since its trace is non-zero), so by adding it
to the algebra we just get $\foc(1;4n) \ltimes \Kee I_0$ from
$\foc(2;4n) \ltimes \Kee I_0$.

\sssec{Summary: The types of Lie superalgebras preserving
non-degene\-ra\-te symmetric forms} In addition to the isomorphisms
$\foo_{\Pi I}(a|b)\simeq\foo_{I\Pi }(b|a)$, there is the only
``occasional" isomorphism intermixing the types of Lie superalgebras
preserving non-degene\-ra\-te symmetric forms:
$\foo_{\Pi\Pi}'(6|2)\simeq\fpe'(4)$.

Let $\widehat{\fg}:=\fg\ltimes\Kee I_0$.  We have the following
types of non-isomorphic Lie (super)algebras:
\begin{equation}\label{oandoo}\renewcommand{\arraystretch}{1.6}
\begin{tabular}{|l|l|}
\hline no relative has Cartan matrix&with Cartan matrix\\
\hline $\foo_{II}(2n+1|2m+1),\ \
\foo_{II}(2n+1|2m)$&$\widehat{\fo\fc(i;2n)},\ \ \fo'(2n+1);\
\ \widehat{\fpe\fc(i;k)}$\\
$\foo_{II}(2n|2m), \ \ \foo_{I\Pi}(2n|2m);\ \
\fo_I(2n);$&$\widehat{\foo\fc(i;2n|2m)},
\ \ \foo_{I\Pi}'(2n+1|2m)$\\
\hline\end{tabular}
\end{equation}

\paragraph{On various versions of the orthogonal Lie algebra, and its prolong, for $p=2$} Let
us begin with $\fg=\fh(2n)$ and $\fh(0|m)$, and for $p=0$ for
simplicity. Both these algebras can be realized on generating
functions (in even and odd indeterminates, respectively) with the
well-known brackets. The component of Lie-degree 0 (in the standard
$\Zee$-grading of $\fg$) is spanned by monomials of degree 2 and is
isomorphic to $\fsp(2n)$ for $\fh(2n)$ and $\fo(m)$ for $\fh(0|m)$.
If we forget the parity of the indeterminates for a moment and look
at the basis of $\fg_0$, the only difference between $\fsp(2n)$ and
$\fo(m)$ is in the fact that the generating functions of the basis
elements of $\fsp(2n)$ contain squares of the indeterminates,
whereas the generating functions of the basis elements of $\fo(m)$
do not contain squares.

Revenons \`a nos moutons, i.e., to $p= 2$ and Lie algebras (no
super!). In this case, as A.~Lebedev explained in \cite{Le1, LeP},
the what he denoted by $\fo$ with various sub- and super-scripts
looks more like the good old $\fsp$, whereas both $\fo_I'$ and
$\fo_\Pi'$ are  the true analogs of the usual $\fo$. Indeed: as
modules over themselves and for $p\neq 2$, we have $\fsp(W)=S^2(W)$,
whereas $\fo(V)=E^2(V)$; while for $p=2$, we have $S^2(V)\supset
E^2(V)$.

Consider the Cartan prolongs of the pairs $(V, \fo(V))$ and $(V,
\fo'(V))$ and realize these prolongs by generating functions. We see
that $(V, \fo(V))_{*, N}$ and $(V, \fo'(V))_{*, N}$ resemble
$\fh(2n)$ and $\fh(0|m)$, respectively.

But the second case can be also interpreted as follows: we {\bf
declared} $N_i=1$ for {\bf all} coordinates of the shearing vector
$N$. Observe that here we are talking about the shearing parameter
for generating functions! The shearing parameter in the realization
of the elements of the algebra by vector fields does not demonstrate
this effect, cf. \cite{LeP}.

One can also take an intermediate road: set $N_i=1$ for SOME $i$,
setting $N_i>1$ for the remaining values of $i$. Then $\fg_0$
becomes isomorphic to something in-between $\fo$ and $\fo'$: In
terms of generating functions, we add (divided) squares of those
indeterminates $x_i$ for which $N_i>1$. In particular, if such an
indeterminate is unique, then $\fo'$ is augmented by ONE element
only.

The cases with restrictions on the coordinates of the shearing
vector of the form $N_i=1$ for some $i$  can also be interpreted as
certain analogs of divergence. We will need several of them. There
are two types of Cartan prolongs of the derived orthogonal Lie
algebras $\fo_B^{(1)}$. These prolongs --- ``little" Hamiltonian Lie
algebras, $\fl\fh_I(n;\uN)$ and $\fl\fh_\Pi(2k;\uN)$
--- consist of vector fields
$A=\mathop{\sum}\limits_{1\leq i\leq n} A_i\partial_i$, elements of
the ``full" Lie algebras $\fh_I(n;\uN)$ and $\fh_\Pi(2k;\uN)$,
satisfying the following conditions:
\begin{equation}\label{hsma}
\begin{array}{ll}
\text{for $\fo^{(1)}_I(n)$}:&\partial_i A_i=0\text{~~for all
$i=1,\dots,n$;}\\
\text{for $\fo^{(1)}_\Pi(2k)$}:&\partial_i
A_{k+i}=\partial_{k+i}A_i=0\text{~~ for all $i=1,\dots,k$.}
\end{array}
\end{equation}

There is also $\fs\fl\fh_\Pi(2k)$, the Cartan prolong of the second
derived Lie algebra $\fo_\Pi^{(2)}(2k)$ consisting of
divergence-free elements of $\fl\fh_\Pi(2k;\uN)$. In \cite{ILL}, we
set:
\begin{equation}\label{defh}
\renewcommand{\arraystretch}{1.4}
\begin{array}{ll}
\fh_I(n;\uN):=(\id, \fo_I(n))_{*, \uN};&
\fh_S(n;\uN):=(\id, \fo_S(n))_{*, \uN};\\
\widetilde\fh_I(n):=(\id, \fc(\fo_I^{(1)}(n)))_{*};&
\widetilde\fh_S(n):=(\id, \fc(\fo_S^{(1)}(n)))_{*},
\end{array}\end{equation}
where $\widetilde\fh$ from \cite{ILL} is the same as $\fs\fl\fh$ in
\cite{LeP}. Now, denote by ${\bf F}(\fl\fle(n;\un|n))$ the
subalgebra of ``half-divergence"-free Hamiltonian vector fields, see
subsec. 3.6 of \cite{ILL} and eq. (2.15) of \cite{LeP}:
\begin{equation}\label{defhalfdiv}
(\id,
\fo_\Pi^{(1)}(2n))_{*}:=\{H_f\mid\sum\nfrac{\del^2f}{\del_{q_i}\del_{p_i}}
=0, \ \un=\un_s\}.
\end{equation}

\ssec{Dynkin diagrams} A usual way to represent simple Lie algebras
over $\Cee$ with integer Cartan matrices is via graphs called, in
the finite dimensional case, {\it Dynkin diagrams}. The Cartan
matrices of certain interesting infinite dimensional simple Lie {\it
super}algebras $\fg$ (even over $\Cee$) can be non-symmetrizable or
have entries belonging to the ground field $\Kee$. Still, it is
always possible to assign an analog of the Dynkin diagram to each
(modular) Lie (super)algebra with Cartan matrix,  provided the edges
and nodes of the graph (Dynkin diagram) are rigged with an extra
information. Although these analogs of the Dynkin graphs are not
uniquely recovered from the Cartan matrix (and the other way round),
they give a graphic presentation of the Cartan matrices and help to
observe some hidden symmetries.

Namely, the {\it Dynkin diagram}\index{Dynkin diagram} of a
normalized $n\times n$ Cartan matrix $A$ is a set of $n$ nodes
connected by multiple edges, perhaps endowed with an arrow,
according to the usual rules (\cite{K}) or their modification, most
naturally and unambiguously formulated by Serganova: compare
\cite{FLS} with vague definitions in \cite{WK, FSS}. In what
follows, we recall these rules, and further improve them to fit the
modular case.

\sssec{Nodes} To every simple root there corresponds
\begin{equation}\label{cm1}
\begin{cases}
\text{a node}\  \mcirc\  &\text{if $p(\alpha_{i})= \ev$ and $
A_{ii}=2$},\\
\text{a node}\  \ast \ &\text{if $p(\alpha_{i}) =\ev$ and
$A_{ii}=\od$};\\
\text{a node}\  \mbullet \ &\text{if $p(\alpha_{i}) =\od$ and
$A_{ii}=1$};\\
\text{a node}\  \motimes \ & \text{if $p(\alpha_{i}) =\od$ and $
A_{ii}=0$},\\
\text{a node}\  \odot\  &\text{if $p(\alpha_{i})= \ev$ and $
A_{ii}=\ev$}.\\
\end{cases}
\end{equation}

The Lie algebras $\fsl(2)$ and $\fo(3)'$ with Cartan matrices $(2)$
and $(\od)$, respectively, and the Lie superalgebra $\fosp(1|2)$
with Cartan matrix $(1)$ are simple.

The Lie algebra $\fgl(2)$ with Cartan matrix $(\ev)$ and the Lie
superalgebra $\fgl(2|2)$ with Cartan matrix $(0)$ are solvable of
$\dim 4$ and $\sdim 2|2$, respectively. Their derived algebras are
the {\it Heisenberg algebra} $\fhei(2):=\fhei(2|0)\simeq\fsl(2)$ and
the {\it Heisenberg superalgebra} $\fhei(0|2)\simeq\fsl(1|1)$ of
(super)dimension 3 and $1|2$, respectively.

\parbegin{Remark} {\it A posteriori} (from the classification of
simple Lie superalgebras with Cartan matrix and of polynomial
growth) we find out that {\bf for $p=0$}, the simple root~$\odot$
can only occur if $\fg(A, I)$ grows faster than polynomially. Thanks
to classification again, if $\dim \fg<\infty$, the simple root
$\odot$ can not occur if $p>3$; whereas for $p=3$, the Brown Lie
algebras are examples of $\fg(A)$ with a simple root of type
$\odot$; for $p=2$, such roots are routine.
\end{Remark}

\sssec{Edges} If $p=2$ and $\dim \fg(A)<\infty$, the Cartan matrices
considered are symmetric. If $A_{ij}=a$, where $a\neq 0$ or 1, then
we rig the edge connecting the $i$th and $j$th nodes by a label $a$.

If $p>2$ and $\dim \fg(A)<\infty$, then $A$ is symmetrizable, so let
us symmetrize it, i.e., consider $DA$ for an invertible diagonal
matrix $D$. Then, if $(DA)_{ij}=a$, where $a\neq 0$ or $-1$, we rig
the edge connecting the $i$th and $j$th nodes by a label $a$.

If all off-diagonal entries of $A$ belong to $\Zee/p$ and their
representatives are selected to be non-positive integers, we can
draw the Dynkin diagram as for $p=0$, i.e., connect the $i$th node
with the $j$th one by $\max(|A_{ij}|, |A_{ji}|)$ edges rigged with
an arrow $>$ pointing from the $i$th node to the $j$th if
$|A_{ij}|>|A_{ji}|$ or in the opposite direction if
$|A_{ij}|<|A_{ji}|$.

\sssec{Reflections} Let $R^+$ be a system of positive roots of Lie
superalgebra $\fg$, and let $B=\{\sigma_1,\dots,\sigma_n\}$ be the
corresponding system of simple roots with some corresponding pair
$(A=A_B,I=I_B)$. Then the set
$(R^+\backslash\{\sigma_k\})\coprod\{-\sigma_k\}$ is a system of
positive roots for any $k\in \{1, \dots, n\}$. This operation is
called {\it the reflection in $\sigma_k$}; it changes the system of
simple roots by the formulas
\begin{equation}
\label{oddrefl}
r_{\sigma_k}(\sigma_{j})= \begin{cases}{-\sigma_j}&\text{if~}k=j,\\
\sigma_j+B_{kj}\sigma_k&\text{if~}k\neq j,\end{cases}\end{equation}
where
\begin{equation}
\label{Boddrefl}B_{kj}=\begin{cases}
-\displaystyle\frac{2A_{kj}}{A_{kk}}& \text{~if~}p_k=\ev,\ \
A_{kk}\neq 0,\text{~and~}
-\displaystyle\frac{2A_{kj}}{A_{kk}}\in \Zee/p\Zee,\\
p-1&\text{~if~}p_k=\ev,\ \  A_{kk}\neq 0\text{~and~}
 -\displaystyle\frac{2A_{kj}}{A_{kk}}\not\in \Zee/p\Zee,\\
-\displaystyle\frac{A_{kj}}{A_{kk}}& \text{~if~}p_k=\od,\ \
A_{kk}\neq 0,\text{~and~}
-\displaystyle\frac{A_{kj}}{A_{kk}}\in \Zee/p\Zee,\\
p-1&\text{~if~}p_k=\od,\ \  A_{kk}\neq 0,
\text{~and~} -\displaystyle\frac{A_{kj}}{A_{kk}}\not\in \Zee/p\Zee,\\
1&\text{~if~}p_k=\od,\ \  A_{kk}=0,\ \ A_{kj}\neq 0,\\
0&\text{~if~}p_k=\od,\ \  A_{kk}=A_{kj}=0,\\
p-1&\text{~if~}p_k=\ev,\ \  A_{kk}=\ev,\ \ A_{kj}\neq 0,\\
0&\text{~if~}p_k=\ev,\ \  A_{kk}=\ev,\ \ A_{kj}=0,\end{cases}
\end{equation}
where we consider $\Zee/p\Zee$ as a subfield of $\Kee$.

The name ``reflection'' is used because in the case of (semi)simple
finite-dimensional Lie algebras this action extended on the whole
$R$ by linearity is a map from $R$ to $R$, and it does not depend on
$R^+$, only on $\sigma_k$. This map is usually denoted by
$r_{\sigma_k}$ or just $r_{k}$. The map $r_{\sigma_i}$ extended to
the $\Ree$-span of $R$ is reflection in the hyperplane orthogonal to
$\sigma_i$ relative the bilinear form dual to the Killing form.

The reflections in the even (odd) roots are said to be {\it even}
({\it odd}).\index{Reflection! odd}\index{Reflection! even!
non-isotropic} \index{Reflection! even! isotropic} A simple root,
and reflection in it, is called {\it isotropic}, if the
corresponding row of the Cartan matrix has zero on the diagonal, and
{\it non-isotropic} otherwise.

If there are isotropic simple roots, the reflections $r_\alpha$ do
not, as a rule, generate a version of the {\it Weyl group} because
the product of two reflections in nodes not connected by one
(perhaps, multiple) edge is not defined.\footnote{The ideas of the
paper \cite{SV} might be helpful here.} These reflections just
connect pair of ``neighboring'' systems of simple roots and there is
no reason to expect that we can multiply two distinct such
reflections. In the general case (of Lie superalgebras and $p>0$),
the action of a given isotropic reflections (\ref{oddrefl}) can not,
generally, be extended to a linear map $R\tto R$. For Lie
superalgebras over $\Cee$, one can extend the action of reflections
by linearity to the root lattice but this extension preserves the
root system only for $\fsl(m|n)$ and $\fosp(2m+1|2n)$, cf.
\cite{Se}.

If $\sigma_i$ is an odd isotropic root, then the corresponding
reflection sends one set of Chevalley generators into a new one:
\begin{equation}
\label{oddrefx} \widetilde X_{i}^{\pm}=X_{i}^{\mp};\ \  \widetilde
X_{j}^{\pm}=\begin{cases}[X_{i}^{\pm},
X_{j}^{\pm}]&\text{if $A_{ij}\neq 0, \ev$},\\
X_{j}^{\pm}&\text{otherwise}.\end{cases}
\end{equation}

\parbegin{Remark} The description of the numbers $B_{ik}$ is
empirical and based on classification \cite{BGL1}: For
infinite-dimensional Lie (super)algebras these numbers might be
different. In principle, in the second, fourth and penultimate
cases, the matrix \eqref{Boddrefl} can be equal to $kp-1$ for any
$k\in\Nee$, and in the last case any element of $\Kee$ may occur.
For $\dim\fg<\infty$, this does do not happen (and it is of interest
to investigate at least the simplest infinite dimensional case ---
the modular analog of \cite{CCLL}).

The values $-\displaystyle\frac{2A_{kj}}{A_{kk}}$ and
$-\displaystyle\frac{A_{kj}}{A_{kk}}$ are elements of $\Kee$, while
the roots are elements of a vector space over $\Ree$. Therefore {\bf
the expressions in the first and third cases in} (\ref{Boddrefl})
{\bf should be understood as} ``{\bf the minimal non-negative
integer congruent to $-\displaystyle\frac{2A_{kj}}{A_{kk}}$ or
$-\displaystyle\frac{A_{kj}}{A_{kk}}$, respectively''. (If
$\dim\fg<\infty$, these expressions are always congruent to
integers.)}

{\bf There is known just one exception: If $p=2$ and
$A_{kk}=A_{jk}$, then} $-\displaystyle\frac{2A_{jk}}{A_{kk}}$ {\bf
should be understood as $2$, not 0.}
\end{Remark}

\paragraph{On neighboring root systems} Serganova \cite{Se} proved (for
$p=0$) that there is always a chain of reflections connecting $B_1$
with some system of simple roots $B'_2$ equivalent to $B_2$ in the
sense of definition \ref{EqSSR}. Here is the modular version of
Serganova's Lemma. Observe that Serganova's statement is not weaker:
Serganova used only odd reflections.

\begin{Lemma}[\cite{LeD}]\label{serg} For any two systems of simple roots
$B_1$ and $B_2$ of any simple finite dimensional Lie superalgebra
with Cartan matrix, there is always a chain of reflections
connecting $B_1$ with $B_2$.\end{Lemma}

\ssec{The Lie (super)algebras of the form $\fg(A)$. Their simple
subquotients $\fg'(A)/\fc$}\label{Ssteps}


\sssec{Over $\Cee$} 
Kaplansky was the first (see his newsletters in \cite{Kapp}) to
discover the exceptional algebras $\fag(2)$ and $\fab(3)$ (he dubbed
them $\Gamma_2$ and $\Gamma_3$, respectively) and a parametric
family $\fosp(4|2; \alpha)$ (he dubbed it $\Gamma(A, B, C))$); our
notation reflect the fact that $\fag(2)_\ev=\fsl(2)\oplus\fg(2)$ and
$\fab(3)_\ev=\fsl(2)\oplus\fo(7)$ ($\fo(7)$ is $B_3$ in Cartan's
nomenclature). Kaplansky's description (irrelevant to us at the
moment except for the fact that $A$, $B$ and $C$ are on equal
footing) of what we now identify as $\fosp(4|2; a)$, a parametric
family of deforms of $\fosp(4|2)$, made
 an $S_3$-symmetry of the parameter manifest (to A.\ A.~Kirillov,
 and he informed us, in 1976).
Indeed, since $A+B+C=0$, and $a\in \Cee\cup\infty$ is the ratio of
the two parameters  remaining after $A$, $B$ and $C$  were
constrained, we get an $S_3$-action on the plane $A+B+C=0$ which in
terms of $a$ is generated by the transformations:
\begin{equation}\label{osp42symm}
a\longmapsto -1-a, \qquad a\longmapsto \frac{1}{a}.
\end{equation}
The other transformations generated by (\ref{osp42symm}) are
\[a\longmapsto -\frac{1+a}{a},\quad a\longmapsto
-\frac{1}{a+1},\quad a\longmapsto-\frac{a}{a+1}.
\]
This symmetry should have immediately sprang to mind since
$\fosp(4|2; a)$ is strikingly similar to $\fwk(3; a)$ found 5 years
earlier, cf. (\ref{wkiso}), and since $S_3\simeq \SL(2; \Zee/2)$.

\sssec{Modular Lie algebras and Lie superalgebras}

\paragraph{$p=2$, Lie algebras} Weisfeiler and Kac \cite{WK} discovered two
new parametric families that we denote $\fwk(3;a)$ and $\fwk(4;a)$
for {\it Weisfeiler and Kac algebras}.\index{$\fwk(3;a)$, Weisfeiler
and Kac algebra} \index{$\fwk(3;a)$, Weisfeiler and Kac
algebra}\index{Weisfeiler and Kac algebra}

\underline{$\fwk(3;a)$, where $a\neq 0, 1$}, of dim 18 is a
non-super version of $\fosp(4|2; a)$ (although no $\fosp$ exists for
$p=2$); the dimension of its simple subquotient $\fwk(3;a)'/\fc$ is
equal to 16; the inequivalent Cartan matrices are:
\[
1)\  \begin{pmatrix} \ev &a &0\\
a&\overline{0}&1\\0&1&\overline{0} \end{pmatrix} ,\
2)\  \begin{pmatrix} \ev &1+a &a\\
1+a&\overline{0}&1\\
a&1&\overline{0} \end{pmatrix}
\]

\underline{$\fwk(4;a)$, where $a\neq 0, 1$}, of $\dim=34$; the
inequivalent Cartan matrices are:
\[
1)\  \begin{pmatrix} \ev &a &0&0\\
a &\overline{0}&1&0\\
0&1&\overline{0}&1\\
0&0&1&\overline{0} \end{pmatrix} ,\
2)\  \begin{pmatrix} \ev &1 &1+a&0\\
1 &\overline{0}& a & 0\\
a+1& a &\overline{0}&a\\
0&0&a&\overline{0} \end{pmatrix}
 ,\
3)\ \begin{pmatrix} \ev &a & 0 &0\\
a &\overline{0}& a+1 & 0\\
0& a+1 &\overline{0}&1\\
0&0&1&\overline{0} \end{pmatrix}
\]

Weisfeiler and Kac investigated also which of these algebras are
isomorphic and the answer is as follows:
\begin{equation}\label{wkiso}
\renewcommand{\arraystretch}{1.4}\begin{array}{l}
\fwk(3;a)\simeq \fwk(3;a')\Longleftrightarrow
a'=\displaystyle\frac{\alpha a+\beta}{\gamma a+\delta},\text{where
$\begin{pmatrix}\alpha&\beta\\ \gamma &\delta\end{pmatrix}\in \SL(2; \Zee/2)$}\\
\fwk(4;a)\simeq \fwk(4;a')\Longleftrightarrow
a'=\displaystyle\frac{1}{a}.
\end{array}
\end{equation}

\paragraph{$p=2$, Lie superalgebras} The same Cartan matrices as for $\fwk$
algebras but with arbitrary distribution of 0's on the main diagonal
correspond to Lie superalgebras $\fbgl(3;a)$ and $\fbgl(4;a)$
discovered in \cite{BGL1}. The conditions when they are isomorphic
are the same as in \eqref{wkiso}, they have the same inequivalent
Cartan matrices, and are considered also only if $a\neq 0, 1$ (since
otherwise they are not simple). We have $\sdim \fbgl(3;a)=10/8|8$
and $\sdim \fbgl(4;a)=18|16$.

\ssec{Yamaguchi's theorem} Let $\fs:=\mathop{\oplus}\limits_{i\geq
-d}\fs_i$ be a simple finite dimensional Lie algebra over $\Cee$.
Let $(\fs_-)_*=(\fs_-, \fg_0)_*$ be the Cartan prolong with the
maximal possible $\fg_0:=\fder_0(\fs_-)$. As is now well-known
\cite{K},
\begin{equation}\label{Zgra}
\begin{minipage}[l]{13cm}
any $\Zee$-grading of the finite dimensional Lie algebra $\fg$
with Cartan matrix is given by a vector $r=(r_1,\dots, r_{\rk\fg})$,
where $r_i\in\Zee$ for all $i$, by setting $\deg X_i^\pm=\pm r_i$.
\end{minipage}
\end{equation}
We say that a grading is \emph{simplest} if $r_{i}=\delta_{ii_0}$
for some $i_0$. The indices $i$ for which $r_{i}=1$ will be called
``selected" (assuming $r_{j}=0$  for all non-selected indices $j$).

\begin{Theorem}[\cite{Y}] Over $\Cee$, equality $(\fs_-)_*=\fs$ holds
almost always. The exceptions (cases where
$\fs=\mathop{\oplus}\limits_{i\geq -d}\fs_i$ is a partial prolong
in $(\fs_-)_*=(\fs_-, \fg_0)_*$) are

{\em 1)} $\fs$ with the grading of depth $d=1$ (in which case
$(\fs_-)_*=\fvect(\fs_-^*)$);

{\em 2)} $\fs$ with the grading of depth $d=2$ and $\dim\fs_{-2}=1$,
i.e., with the ``contact'' grading, in which case
$(\fs_-)_*=\fk(\fs_-^*)$ (these cases correspond to ``selecting" the
nodes on the Dynkin graph connected with the node representing the
maximal root on the extended graph);

{\em 3)} $\fs$ is either $\fsl(n+1)$ or $\fsp(2n)$ with the grading
determined by ``selecting'' the first and the $i$th of simple
coroots, where $1<i<n$ for $\fsl(n+1)$ and $i=n$ for $\fsp(2n)$.
(Observe that $d=2$ with $\dim \fs_{-2}>1$ for $\fsl(n+1)$ and $d=3$
for $\fsp(2n)$.)

Moreover, the equality $(\fs_-, \fs_0)_*=\fs$ also holds almost
always. The cases where the equality fails (the ones where a
projective action is possible) are $\fsl(n+1)$ or $\fsp(2n+2)$ with
the grading determined by ``selecting'' only one (the first) simple
coroot; $\fs=\fvect(n)$ or $\fk(2n+1)$, respectively.
\end{Theorem}
\sssbegin{Remark} First, Yamaguchi's cases (for $p=0$) where the CTS
prolongs return the initial algebra are precisely the cases where
restrictions on $\un$ are imposed if we pass to $p>0$. More exactly,
to describe the \emph{complete} prolong, NO restrictions should be
imposed on $\underline{N} $.

It so happens that (even if $p=0$) certain indeterminates, in terms
of which the CTS prolong is described, can not enter in degrees
greater than something. For $p>0$, this imposes certain restrictions
on $\un$ dictated by the very structure of the Lie algebra whose
non-positive part we are prolonging.

For example, if we write $\underline N=(1,n,1)$, this does not mean
that WE have imposed any constraints on the first and third
coordinates of $\underline N$; it is the non-positive (or negative,
depending on the problem) part of the algebra to be prolonged
imposes these constraints on these coordinates.

Therefore, not only in the case where there are no restrictions on
$\underline N$ but also in all cases where at least one of the
coordinates of $\underline N$ is not restricted, the COMPLETE
prolong is of infinite dimension. The space $\fg_1$ can not generate
the complete prolong, or any part of it with sufficiently great
value of at least one of coordinates of $\underline N$. The algebra
generated by $\fg_1$ gives us restrictions on coordinates of
$\underline N$ imposed by $\fg_1$.
\end{Remark}

\section{Simple Lie algebras as CTS-prolongs of the non-positive parts of $\fg(A)$}
For the definition of the grading vector $r$ in tables below, see
\eqref{Zgra}; for a $\Zee$-graded Lie algebra $\fg=\oplus \fg_i$, we
set $\fg_{\leq0}:=\mathop{\oplus}\limits_{i\leq 0}\fg_i$. Consider
Cartan matrices as their size grows. In tables \eqref{rank3} and
\eqref{rank4}, we provisionally (until we identify the algebra
$\fg_{*, \un}$ with a known algebra)  denote the prolongs $\fg_{*,
\un}$ with $\dim\fg_-=D$ by $D(D;\un)$.

\ssec{Size=1}
\begin{equation}\label{rank1}
\renewcommand{\arraystretch}{1.4}
\begin{tabular}{|c|c|c|l|}
\hline $\fg$&Cartan matrix&$r$&Prolong of $\fg_{\leq0}$ for this $r$\\
\hline
 $\begin{matrix}\fo_{\Pi}'(3)
\end{matrix}$&$ \begin{pmatrix}
    \od
  \end{pmatrix}$&$1$&
$\fvect(1;\un)$\\
\hline
\end{tabular}
\end{equation}

\normalsize

\ssec{Size=2}
 \footnotesize

\begin{equation}\label{rank2}
\renewcommand{\arraystretch}{1.4}
\begin{tabular}{|c|c|c|l|}
\hline $\fg$&Cartan matrix&$r$&Prolong of $\fg_{\leq0}$ for this
$r$\\ \hline $\begin{matrix}\fsl(3)
\end{matrix}$&$ \begin{pmatrix}
    \ev &1 \\
    1& \ev
  \end{pmatrix}$&$\begin{matrix}(10)\\
(01)\end{matrix}$& $\begin{matrix}\fvect(2;\un)\\
\fvect(2;\un)\end{matrix}$\\[3mm] \hline $\begin{matrix}\fo_{\Pi}'(5)
\end{matrix}$&$ \begin{pmatrix}
    \ev & 1 \\
    1 & \overline{1}
  \end{pmatrix}$&$\begin{matrix}(10)\\
(01)\end{matrix}$&
$\begin{matrix}\mathfrak{B}(3;\un)=\fh_\Pi(3;\un)\\
\mathfrak{C}(3;\un)=\fk(3;\un)\end{matrix}$\\
\hline
\end{tabular}
\end{equation}

\normalsize

\paragraph{What $\mathfrak{B}(2m-1;\un)$ and
$\mathfrak{C}(2m-1;\un)$ are} For $m=2$, consider the generating
functions in $x_1,x_2,x_3$ with $\un=(1,n,1)$ and the Poisson
bracket:
$$
\{f,g\}=\frac{\partial f}{\partial x_1}\frac{\partial g}{\partial
x_3}+ \frac{\partial f}{\partial x_3}\frac{\partial g}{\partial x_1}
+ \frac{\partial f}{\partial x_2}\frac{\partial g}{\partial x_2}.
$$
Having factorized modulo center (generated by constants) we get
$\fh_\Pi(3;\un)$, and its elements can be represented as
$$
H_f=\frac{\partial f}{\partial x_1}\partial_{x_3}+ \frac{\partial
f}{\partial x_3}\partial_{x_1} + \frac{\partial f}{\partial
x_2}\partial_{x_2}.
$$

Setting $\deg_{Lie}=\deg-2$, where $\deg x_i=1$ for all $i$, we see
that the height of $\fh_\Pi(3;\un)$ relative the grading
$\deg_{Lie}$ is equal to $2^n-1$, and
\[
\fg_k=\Span(x_2^{(k+2)}, \quad x_2^{(k+1)}x_1, \quad x_2^{(k+1)}x_3,
\quad x_2^{(k)}x_1x_3)\ \text{for $0\le k<2^n-1$,}
\]
whereas  $\dim \fg_{2^n-1}=1$. Since the same arguments hold for
$m>2$ as well, we arrive at the following verdict:
\[
\mathfrak{B}(2m-1;\un)\cong \fh_\Pi(2m-1;\un).
\]
By comparing non-positive parts of the $\Zee$-graded Lie algebras
$\mathfrak{C}(2m-1;\un)$ we deduce:
\[
\begin{array}{l}\mathfrak{C}(3;\un)\cong \fk(3;\un),\\
\mathfrak{C}(2m-1;\un)=\fo_\Pi(2m+1)\text{~~for $m>2$.}
\end{array}
\]

\ssec{Size=3} \sssec{Derivarions and central extensions} When the
size of Cartan matrix is $\geq  3$ it might be non-invertible, and
hence the algebras $\fg(A)$ have to be replaced, see the left column
of table \eqref{rank3}, with their quotients modulo center; the same
applies to the derived of $\fg(A)$. These latter algebras may have,
for $A$ of small rank, non-trivial outer derivations and in order to
identify the ``extra" elements of CTS-prolongs with {\bf some} of
these derivations we have to know all of them. For example, we
already know that $\fpsl(4)$ is a desuperization of
$\fpsl(2|2)\simeq\fh'(0|4)$, and hence has at least three outer
derivations of degrees $\pm 4$ and $0$ with respect to the grading
defined on Chevalley generators by setting $\deg X_i^\pm=\pm 1$
(these cocycles turn $\fh'(0|4)$ into, respectively, $\fh(0|4)$
--- twice, and $\fpgl(2|2)\simeq \fh'(0|4)\subplus \Kee E$,
where $E=\sum \theta_i\partial_i$ is the Euler operator). In reality
there are 7 central extensions, see \cite{BGLL2}.

For reasons explained in \cite{BGLL2}, there should be no less than
7 outer derivations of $\fpsl(4)$; there are precisely 7 of them,
see \cite{BGLL2}.

\sssec{The table} Explanations of isomorphisms in the right-most
column of table \eqref{rank3} are given under it. For a description
of the Buttin algebras $\fb_{\lambda}(m)$ and their various
Weisfeiler regradings $\fb_\lambda(m;r)$ over $\Cee$, see
\cite{LSh}; their analogs for $p=2$ are described in \cite{LeP}. In
particular, $\fb_\lambda(2;2; \un)$ is the nonstandard regrading of
$\fb_\lambda(2; \un)$ corresponding to $\deg \xi_i=0$ for both odd
indeterminates $\xi_i$, see \cite{LSh}. For $p=2$, the odd
indeterminates correspond to those of the desuperization, whose
degree can not exceed 1.

 \tiny

\begin{equation}\label{rank3}
\renewcommand{\arraystretch}{1.4}
\begin{tabular}{l}
\begin{tabular}{|c|c|c|l|}
\hline $\fg$&Matrix $A$&$r$&Prolong of $\fg_{\leq0}$ for this $r$\\
\hline $\begin{matrix}\fpgl(4)
\end{matrix}$&$ \begin{pmatrix}
    \ev &1&0 \\
    1& \ev&1\\
    0& 1&\ev
  \end{pmatrix}$&$\begin{matrix}(100)\\
(010)\\
(001)\end{matrix}$& $\begin{array}{l}\fvect(3;\un)\\
{\mathbf F}(\fh(0|4)\ltimes \Kee E), \text{~~where $E=\sum\theta_i\partial_{\theta_i}$}\\
\fvect(3;\un)\end{array}$\\[3mm]
\hline $\begin{matrix}\fpsl(4)
\end{matrix}$&$ \begin{pmatrix}
    \ev &1&0 \\
    1& \ev&1\\
    0& 1&\ev
  \end{pmatrix}$&$\begin{matrix}(100)\\
(010)\\
(001)\end{matrix}$& $\begin{matrix}\fsvect(3;\un)\\
\fpgl(4)\\
\fsvect(3;\un)\end{matrix}$\\[3mm] \hline $\begin{matrix}\fo'_{\Pi}(7)
\end{matrix}$&$ \begin{pmatrix}
    \ev & 1& 0 \\
    1 & \ev& 1\\
    0 & 1 & \od
  \end{pmatrix}$&$\begin{matrix}(100)\\
(010)\\
(001)\end{matrix}$& $\begin{matrix}
5(5;N)=\fh_I(5;1,1,n,1,1), \text{ see \cite{ILL}, Cor. 3.9.2}\\
\fo'_{\Pi}(7)\ltimes 2\  \mathrm{outer\  derivations}\\
\fo'_{\Pi}(7)\ltimes 3\  \mathrm{outer\  derivations}\\
\end{matrix}$\\
\hline $\begin{matrix}\fwk'(3;\lambda)/\fc
\end{matrix}$&$ \begin{matrix}\begin{pmatrix}
    \ev & \lambda& 0 \\
    \lambda & \overline{0}& 1\\
    0 & 1 & \ev
  \end{pmatrix}\\
\begin{pmatrix}
    \ev & 1+\lambda& \lambda \\
    1+\lambda & \overline{0}& 1\\
    \lambda & 1 & \ev
  \end{pmatrix}
  \end{matrix}$&$\begin{matrix}(100)\\
(010)\\
(001)\\
(100)\\
(010)\\
(001)\end{matrix}$& $\begin{array}{l}
\mathfrak{G}(5;\lambda;\underline{N}
)={\bf F}(\fb_{\lambda}(2;2;
\un))\\
\fwk'(3;\lambda)/\fc\\
\mathfrak{G}(5;\lambda^{-1};\un)={\bf F}(\fb_{\lambda^{-1}}(2;2;
\un))\\
\widetilde{\mathfrak{G}}(5;\lambda;\un)={\bf F}(\fb_\lambda(2;2;
\un))\\
\widetilde{\mathfrak{G}}(5;\lambda^{-1};\un)={\bf
F}(\fb_{\lambda^{-1}}(2;2;
\un))\\
\widetilde{\mathfrak{G}}(5;\frac{\lambda}{1+\lambda};\un)={\bf
F}(\fb_{{\lambda}/(1+\lambda)}(2;2;
\un))\end{array}$\\
\hline \hline
\end{tabular}\\
For $\fwk(3;\lambda)/\fc$, the prolong of the non-positive part is
the same as for the respective line for $\fwk'(3;\lambda)/\fc$\\
plus one outer derivation, except for the second line of the first
Cartan matrix in which case\\
it is isomorphic to ${\bf F}(\fk(1;\underline{1}|4))$ independently
of parameter $\lambda$.\\\hline\end{tabular}\end{equation}

\normalsize

\sssec{Elucidating table \eqref{rank3}: realization of $5(5;N)$ by
vector fields} We have $\fg_{0}\simeq \fo'(5)\oplus\fc$, where
$\fc=\Kee z$ and $\fg_{-1}$ is the tautological $\fg_{0}$-module.
This case was studied in \cite{ILL}. Recall that the shearing
parameter is of the form $\un=(1,1,n,1,1)$.

\paragraph{A description of $\mathfrak{G}(5;\un;\lambda):=(\fg_-,\fg_0)_{*, \uN}$
for $\fg=\fwk(3; \lambda)$ with the first Cartan matrix}

The grading $r=(100)$ gives $\dim(\fwk(3; \lambda)_-)=4$. The
realization by vector fields is as follows: \footnotesize
\begin{equation}\label{fG(5N1)}\tiny
\renewcommand{\arraystretch}{1.4}
\begin{tabular}{|l|l|} \hline
$\fg_{i}$&the generators  \\ \hline $\fg_{-1}$&$\partial_1,\
\partial_2,\ \partial_3,
\  \partial_4$\\
\hline $\fg_{0}\simeq \fsl(3)$& $Y_2=x_1\partial_2+x_3\partial_4,\
Z_2=\lambda x_2\partial_1+(1+\lambda )x_4\partial_3,\
H_2=[Z_2,Y_2],$\\
& $Y_3=x_2\partial_3,\ Z_3=x_3\partial_2,\ H_3=[Z_3,Y_3],\
Y_4=[Y_2,Y_3],\  Z_4=[Z_2,Z_3]$
   \\
\hline
\end{tabular}
\end{equation} \normalsize
Our computation shows that $\un=(1,n,m,1)$.

\paragraph{A description of
$\widetilde{\mathfrak{G}}(5;\un;\lambda):=(\fg_-,\fg_0)_{*, \uN}$
for $\fg=\fwk(3; \lambda)$  with the second Cartan matrix}

The grading $r=(100)$ gives $\dim(\fwk(3; \lambda)_-)=4$. The
realization by vector fields is as follows: \footnotesize
\begin{equation}\label{fGbis(5N2)}\tiny
\renewcommand{\arraystretch}{1.4}
\begin{tabular}{|l|l|} \hline
$\fg_{i}$&the generators  \\ \hline $\fg_{-1}$&$ \partial_1,\
 \partial_2,\
 \partial_3, \   \partial_4$\\
\hline $\fg_{0}\simeq \fsl(3)$& $Y_2=x_1\partial_2+x_3\partial_4,\
Z_2=(\lambda +1)\, x_2\partial_1+\lambda \, x_4\partial_3,\
H_2=[Z_2,Y_2],\ Y_3=x_1\partial_3+\frac{\lambda
   +1}{\lambda }\, x_2\partial_4$\\
& $Z_3=\lambda(x_3\partial_1+x_4\partial_2), \  Y_4=[Y_2,Y_3], \
Z_4=[Z_2,Z_3]$
   \\
\hline
\end{tabular}
\end{equation} \normalsize Our computation
shows that $\un=(n,1,1,m)$.

\paragraph{Further elucidating Table \eqref{rank3}} There is no mistake in
the description of prolongs for $\fo'(7)$. The gradings $(010)$ and
$(001)$ give the algebra plus different number of linearly
independent outer derivations. There is no reason why prolongs of
negative parts relative different gradings must yield the same
number of outer derivations.

Let us explain where 2 or 3 outer derivations in line 2 (concerning
$\fo'_\Pi(7)$) come  from. The heart of the matter lies in
generating functions $p_i^2$ and $q_i^2$ one disregards when
considering $\fo'$ instead of $\fo$. The generating functions which,
in the regraded, as $r$ varies, algebra might appear in the
components of non-positive degree can not, of course, appear. But
the ones that should appear in components of positive degree may
appear and do so.

This is most clear if you consider the weights for $p=0$. Obviously,
for $p=2$ several distinct weights coincide (modulo 2), but for us
it is only important that we can realize the initial algebra to be
regraded on generating functions (homogeneous degree 2 polynomials
in 7 indeterminates).

The weights (roots) of the elements in the initial algebra are
well-known: $\pm \eps_i\pm\eps_j$ and $\pm \eps_i$. Accordingly, the
simple roots are $$ \alpha_1=\eps_1-\eps_2, \quad
\alpha_2=\eps_2-\eps_3, \quad\alpha_3=\eps_3.$$

Setting $\deg X_{\alpha_1}=1$ whereas $\deg X_{\alpha_i}=0$ for
$i=2,3$, we get $\fg_0=\fo'_\Pi(5)\oplus \Kee$, the depth is equal
to 1, and the prolong is isomorphic to $\fh$.

Setting $\deg X_{\alpha_2}=1$ whereas $\deg X_{\alpha_i}=0$ for
$i\neq 2$, we get the algebra of depth 2, and the component $\fg_1$
contains elements of weight $\eps_1$ and $\eps_2$. If now we set
\[
\text{weight$(p_i)=1$ and hence $\text{weight}(q_i)=-1$,}
\]
then $p_1^2$ and $p_2^2$ belong to $\fg_2$ and personify 2 outer
derivations of $\fg$ (since, e.g., $[p_1^2,q_1\theta]=p_1\theta\in
\fg_1$). In this grading the element of weight $\eps_3$ lies in
$\fg_0$, and so $p_3^2$ should also lie in $\fg_0$, which is
impossible by hypothesis.

Now, setting  $X_{\alpha_3}=1$ whereas $\deg X_{\alpha_i}=0$ for
$i\neq 3$, the component $\fg_1$ contains all the three elements of
weight $\eps_i$. Accordingly, all the three $p_i^2$ lie in $\fg_2$
(and serve as outer derivations).

\paragraph{$(\fwk(\lambda;3)/\fc)_{\leq 0}\simeq{\bf
F}(\fk(1;\underline{1}|4)))_{\leq 0}$ for $r=(010)$ and the first
Cartan matrix} This case is an exceptional one and does not fit the
pattern of the other gradings because in this case the isomorphism
of non-positive parts of two algebras occurs and therefore their
prolongs coincide. We have
\begin{equation}\label{fG(5N)}\tiny
\renewcommand{\arraystretch}{1.4}
\begin{tabular}{|l|l|} \hline
$\fg_{i}$&the generators  \\ \hline
$\fg_{-2}$&$w_{1}=\partial_1$\\\hline $\fg_{-1}$&$
w_{2}=\partial_2,\ w_{3}=\partial_3,
\; w_{4}=x_3\partial_1+\partial_4,\ w_5=x_2\partial_1+\partial_5$ \\
\hline $\fg_{0}\simeq \fo'_\Pi(4)\oplus \fc\simeq$& $X_1^-=
x_2\partial_3+ x_4\partial_5,\ X_1^+=\lambda
(x_3\partial_2+x_5\partial_4),$\\
& $X_2^-=x_2\, x_3\partial_1+x_2\partial_4+x_3\partial_5,\
X_2^+=x_4\, x_5\partial_1+x_4\partial_2 +x_5\partial_3$\\
$\fhei(4)\ltimes 2 \text{ outer}$& $H_1=[X_1^+, X_1^-]=\lambda H_2,\
H_3=(\lambda +1)\, x_1\partial_1 + \lambda \,
   x_3\partial_3 + x_4\partial_4+(\lambda +1)\, x_5\partial_5$,\\
 derivatives& $H_2=[X_2^+, X_2^-]=x_2\partial_2+ x_3\partial_3 + x_4\partial_4+x_5\partial_5$,\
   $d=x_1\partial_1+x_3\partial_3+x_5\partial_5$\\
\hline
\end{tabular}
\end{equation}
An isomorphism between the non-positive parts of ${\bf
F}(\fk(1;\underline{1}|4)))$, see \cite{LeP}, and
$\fwk(\lambda;3)/\fc$ goes as follows, where we briefly write $f$
instead of contact vector field $K_f$:
\[
\begin{array}{lll}\begin{array}{l}
w_2 \longleftrightarrow
\xi_1\\
w_3 \longleftrightarrow
 \xi_2\\
w_4 \longleftrightarrow \eta_2\\
w_5 \longleftrightarrow \eta_1\\
\end{array}&\begin{array}{l}
X_2^- \longleftrightarrow \xi_2\eta_1\\
X_2^+ \longleftrightarrow \lambda \xi_1\eta_2\\
X_1^- \longleftrightarrow \eta_1\eta_2\\
X_1^+ \longleftrightarrow \xi_1\xi_2\\
\end{array}&\begin{array}{ll}
H_1 &\longleftrightarrow \lambda (\xi_1\eta_1+ \xi_2\eta_2)\\
H_3 &\longleftrightarrow
 \lambda\xi_2\eta_2+(1+\lambda)t\\
H_2 &\longleftrightarrow \xi_1\eta_1+\xi_2\eta_2\\
d &\longleftrightarrow \xi_2\eta_2+t\end{array}\end{array}
\]
Recall that the contact bracket corresponding to the bracket $[K_f,
K_g]$ is
\begin{equation}\label{Konbr}
\{f,g\}_{k.b.}=\pderf{f}{t}(1-E')(g)+(1-E')(f)\pderf{g}{t}+\{f,g\}_{P.b.},
\end{equation}
where $E'=\sum \xi_i\partial_{\xi_i}$ and the Poisson bracket is:
\[\{f,g\}_{P.b.}=\pderf{f}{\xi_1}
\pderf{g}{\eta_1}+\pderf{g}{\xi_1}
\pderf{f}{\eta_1}+\pderf{f}{\xi_2}
\pderf{g}{\eta_2}+\pderf{g}{\xi_2}\pderf{f}{\eta_2}.\] The
non-positive parts determine the isomorphic prolongs that do not
depend on $\lambda$; the algebra $\fwk(3; \lambda)$ is a subalgebra
of ${\bf F}(\fk(1;\underline{1}|4)))$, a partial prolong.

\paragraph{What $\mathfrak{G}(5;\un;\lambda)$ and
$\widetilde{\mathfrak{G}}(5;\un;\lambda)$ are isomorphic to} For
various Cartan matrices and various simplest regradings $r$, we have
$\dim(\fwk(3; \lambda)_-)=4$, see \eqref{rank3}. The realization by
vector fields is  given in eqs. \eqref{fG(5N1)} and
\eqref{fGbis(5N2)}. Let us compare the representations of $\fg_0$ in
$\fg_{-1}$ in these two cases.

Let us begin with $\mathfrak{G}(5;\un;\lambda)$. The highest weight
vector in $\fg_{-1}$ is the one killed by all the $Y_i$. From eq.
\eqref{fG(5N1)} we see that this is $\partial_4$.

Now look how the $Z_i$ act on it:
\begin{equation}\label{line}
\begin{CD}
\partial_4 @> Z_2 >> (1+\lambda)\partial_3 @>Z_3 >> (1+\lambda)\partial_2 @> Z_2 >>
\lambda(\lambda+1)\partial_1 @> Z_2,\ Z_3 >> 0\\
@VVZ_3V     @VVZ_2V         @VVZ_3V         @.\\
0 @. 0 @. 0 @.
\end{CD}
\end{equation}

To compute the weights, we need an explicit form of the $H_i$:
\[
H_2=(\lambda+1)(x_3\partial_3+x_4\partial_4)+
\lambda(x_1\partial_1+x_2\partial_2);\quad
H_3=x_2\partial_2+x_3\partial_3.
\]


Thus, the weight diagram of our representation is as follows:
\begin{equation}\label{tab1}
\begin{tabular}{|l|l|l|l|}
 \hline
$\partial_4$ & $\partial_3$ & $\partial_2$& $\partial_1$\\
\hline $(\lambda+1,0)$&$(\lambda+1,1)$ & $(\lambda,1)$ & $(\lambda,0)$\\
\hline
\end{tabular}
\end{equation}

Now, pass to $\widetilde{\mathfrak{G}}(5;\un;\lambda)$. In order to
avoid a confusion with the previous discussion, let us denote the
basis elements of $\fg_0$ by small letters: $y_i, z_i, h_i$. This
$y_3=\lambda x_1\partial_3+(\lambda+1)x_2\partial_4$ differs from
that of \eqref{fGbis(5N2)} by a factor $\lambda$ and
$z_3=x_3\partial_1+x_4\partial_2$ differs from that of
\eqref{fGbis(5N2)} by a factor $\lambda^{-1}$. Neither $h_3$, nor
Cartan matrix are affected.

We have:
$$
h_2=(\lambda+1)(x_1\partial_1+x_2\partial_2)+\lambda(x_3\partial_3+x_4\partial_4),
\quad
h_3=\lambda(x_1\partial_1+x_3\partial_3)+(\lambda+1)(x_2\partial_2+x_4\partial_4).
$$

Assuming the $y_i$ to be positive root vector we again see that
$\partial_4$ is a highest weight vector. The $z_i$ act on it as
follows:
\begin{equation}\label{square}
\begin{CD}
\partial_4 @> z_2 >> \lambda\partial_3 @>z_3 >> \lambda\partial_1 @> z_2,\ z_3 >> 0\\
@VVz_3V    @.        @.\\
\partial_2 @> z_2 >>  (\lambda+1)\partial_1 @> z_2,\ z_3 >> 0
\end{CD}
\end{equation} and the weight diagram of the representation is of the
following form:
\begin{equation}\label{tab2}
\begin{tabular}{|l|l|l|l|}
 \hline
$\partial_4$ & $\partial_3$ & $\partial_2$& $\partial_1$\\
\hline $(\lambda,\lambda+1)$&$(\lambda,\lambda)$ & $(\lambda+1,\lambda+1)$ & $(\lambda+1,\lambda)$\\
\hline
\end{tabular}
\end{equation}


Now, let us change the basis of $\fg_{-1}$. First, let us
interchange $x_2$ with $x_3$. We see that
\[
h_3=\lambda(x_1\partial_1+x_2\partial_2)+(\lambda+1)(x_3\partial_3+x_4\partial_4)=H_2
\]
whereas
\[
h_2=(\lambda+1)(x_1\partial_1+x_3\partial_3)+\lambda(x_2\partial_2+x_4\partial_4).
\]
Then $h=h_2+h_3=x_1\partial_1+x_4\partial_4$.

Now, let us interchange $x_1$ with $x_2$, and $x_3$ with $x_4$. Then
$h_2$ does not vary, whereas $h$ turns into
$x_2\partial_2+x_3\partial_3=H_3$.

As a result, we have performed the permutation of indeterminates
\begin{equation}\label{perm}\begin{pmatrix}
1&2&3&4\\
 2&4&1&3
\end{pmatrix}\end{equation}which sends the other vectors of
$\mathfrak{G}(5;\un;\lambda)_0$ as follows:
\[
\renewcommand{\arraystretch}{1.4}
\begin{array}{ll}
y_3\mapsto \lambda x_2\partial_1+(\lambda+1)x_4\partial_3=Z_2, &
z_3\mapsto x_1\partial_2+x_3\partial_4=Y_2, \\
y_6\mapsto x_2\partial_3=Y_3,
& z_6\mapsto x_3\partial_2=Z_3,\\
y_2\mapsto x_2\partial_4+x_1\partial_3=Y_5, & z_2\mapsto
(\lambda+1)x_4\partial_2+\lambda x_3\partial_1=Z_5. \end{array}
\]

In order to see how vector diagram \eqref{line} turns into
\eqref{square}, pass to another ``Borel" subalgebra (with positive
generators $Z_2$ and $Y_5$ and respective negative ones $Y_2$ and
$Z_5$) and change basis according to eq. \eqref{perm}. Now
$\partial_3$ becomes a highest weight vector and the images of
$\fg_0$ in $\fgl(\fg_{-1})$ coincide. The Cartan prolongs of these
two pairs $(\fg_{-1}, \fg_0)$ are isomorphic for $\fG$ and
$\widetilde\fG$.

For illustration, let us express the elements of $\fg_0$ (for $\fG$,
not for $\widetilde\fG$) by matrices. The most simple form is
obtained in the basis $\partial_2,
\partial_3,\partial_1,\partial_4$:
\[
\begin{pmatrix}
A&B\\
C&D
\end{pmatrix}, \ \text{ where $A$ and $C$ are arbitrary,}
$$
$$
D=\begin{pmatrix}
   \lambda\cdot\tr A & 0\\
   0 &  (\lambda+1)\cdot\tr A
\end{pmatrix},\
\text{ and if $C=\begin{pmatrix}
            \alpha &\beta \\
            \gamma & \delta
            \end{pmatrix},$ then $B=\begin{pmatrix}
                               \delta  &\nfrac{\lambda+1}\lambda\cdot\beta \\
                                \gamma &\nfrac{\lambda+1}\lambda\cdot\alpha
                                \end{pmatrix}.$}
                                \]
It is not difficult to verify that for
$\fg:=\widetilde{\mathfrak{G}}(5;\un;\lambda)$, the $\fg_0$-action
on $\fg_{-1}$ is precisely the $\fvect(2;\un_s)$-action on the space
of $\lambda$-densities; observe that $\fsl(3)={\bf
F}(\fsl(1|2))\cong{\bf F}(\fvect(0|2))$.

In other words, the prolong is a desuperization of $\fb_\lambda(2;2;
\un)$, the nonstandard regrading of $\fb_\lambda(2; \un)$
corresponding to $\deg \xi_i=0$ for both odd indeterminates $\xi_i$,
see \cite{LSh}. In the version of the prolong we consider, both free
shearing parameters (corresponding to the even indeterminates) are
taken equal to  2. {\bf Here} is the place where the difference
between $\un$ for the generating functions, which correctly
describes the algebra, and $\un$  for coefficients of vector fields,
which yields a wrong description since some information becomes
lost, see ``non-existent generating functions" $\theta^2$ in
\cite{LeP}.

\parbegin{Remarks} 1) This simple Lie algebra --- prolong of $\fwk(3; \lambda)$
--- first appeared  (without interpretation, in components) in
\cite{Bro}. It is difficult not to admire the computational skill of
Brown and compare the difficult calculation made by bare hands in
\cite{Bro} with easiness brought to us by code \emph{SuperLie}.

2) Although the Lie algebra $\fwk(3; \lambda)$ is defined for
$\lambda\neq 0,1$, its prolong $\fb_\lambda(2;2; \un)$ is
well-defined for all values of $\lambda$. For $\lambda= 0,1$, these
prolongs are not simple: for $\lambda= 0$, it has a center the
quotient modulo which is simple, for $\lambda= 1$, it has a simple
ideal of codimension 1.
\end{Remarks}

\ssec{Size=4}

 \tiny

\begin{equation}\label{rank4}
\renewcommand{\arraystretch}{1.6}
\begin{tabular}{|c|c|c|l|}
\hline $\fg$&Cartan matrix&$r$&Prolong of $\fg_{\leq0}$ for this
$r$\\
\hline $\begin{matrix}\mathfrak{oc}(1;8)/\fc\ltimes \Kee I_0
\end{matrix}$&$ \begin{pmatrix}
    \ev & 1& 0 & 0 \\
    1& \ev & 1 & 1 \\
    0& 1& \ev & 0\\
    0 & 1 & 0 &\ev
  \end{pmatrix}$&$\begin{matrix}(1000)\\
(0100)\\
(0010)\\
(0001)\end{matrix}$& $\begin{matrix} 6(6)=\widetilde \fh_\Pi(6)\ltimes \text{2 outer derivatives}\\
9(9;\underline{n, 1, ..., 1})\ltimes 2\, \mathrm{outer\  derivations}\\
6(6)=\widetilde \fh_\Pi(6)\ltimes \text{2 outer derivatives}\\
6(6)=\widetilde \fh_\Pi(6)\ltimes \text{2 outer derivatives}\end{matrix}$\\
[3mm] \hline $\begin{matrix}\mathfrak{oc}(1;8)/\fc =\fo^{(2)}_\Pi(8)/\fc
\end{matrix}$&$ \begin{pmatrix}
    \ev & 1& 0 & 0 \\
    1& \ev & 1 & 1 \\
    0& 1& \ev & 0\\
    0 & 1 & 0 &\ev
  \end{pmatrix}$&$\begin{matrix}(1000)\\
(0100)\\
(0010)\\
(0001)\end{matrix}$& $\begin{matrix}\widetilde{\fh}_\Pi(6), \text{ see \cite{ILL}, eq. (2.6)}\\
9(9;\underline{n, 1, ..., 1})=\fii\fr(9;\un),  \text{ see eq. \eqref{ir}}\\
\widetilde{\fh}_\Pi(6), \text{ see \cite{ILL}, eq. (2.6)}\\
\widetilde{\fh}_\Pi(6), \text{ see \cite{ILL}, eq. (2.6)}\end{matrix}$\\
[3mm] \hline $\begin{matrix}\fsl(5)
\end{matrix}$&$ \begin{pmatrix}
    \ev & 1 & 0& 0 \\
    1 & \ev & 1 & 0 \\
    0 & 1 & \ev & 1 \\
    0 & 0 & 1 & \ev
  \end{pmatrix}$&$\begin{matrix}(1000)\\
(0100)\\(0010)\\(0001)\end{matrix}$&
$\begin{matrix}\mathfrak{vect}(4;\un)\\
\fsl(5)\\
\fsl(5)\\
\mathfrak{vect}(4;\un)
\end{matrix}$\\ \hline
$\mathfrak{o}_{\Pi}'(9)$ & $
\begin{pmatrix}
    \ev & 1 & 0& 0 \\
    1 & \ev & 1 & 0 \\
    0 & 1 & \ev & 1 \\
    0 & 0 & 1 & \od
  \end{pmatrix}$&$\begin{matrix}(1000)\\(0100)\\(0010)\\(0001)\end{matrix}$ &
  $\begin{matrix}\fh_I(7;1,1,1,n,1,1,1), \text{ see \cite{ILL}, Cor. 3.9.2}\\
\mathfrak{o}_{\Pi}'(9)\ltimes 2\  \mathrm{outer\  derivations}\\
\mathfrak{o}_\Pi'(9) \ltimes 3\  \mathrm{outer\  derivations} \\
\mathfrak{o}_\Pi'(9)  \ltimes 4\  \mathrm{outer\  derivations}
\end{matrix}$\\ \hline

$\mathfrak{wk}(4,a)$  &

$ \begin{matrix} \begin{pmatrix}
    \ev & a & 0& 0 \\
    a & \ev & 1 & 0 \\
    0 & 1 & \overline{0} & 1 \\
    0 & 0 & 1 & \ev
  \end{pmatrix}

\\
\text{the other matrices}
  \end{matrix}$&$\begin{matrix}
(1000)\\
(0100)\\
(0010)\\
(0001)\\
\text{any of the above}\end{matrix}$& $\begin{matrix}
\mathfrak{wk}(4;a)\\
11(11;\un;a)={\bf F}(\fm\fb(3;\un|8))\\
\mathfrak{wk}(4;a)\\
\mathfrak{wk}(4;a)\\
\mathfrak{wk}(4;a)
\end{matrix}$\\ \hline
\end{tabular}
\end{equation}
\normalsize

\sssec{The Lie algebra $11(11;\un;a)$ is a desuperization of
$\fm\fb(3;\un|8)$}\label{sssmb} The grading $r=(0100)$ gives
$\dim(\fwk(3; a)_-)=11$. The CTS prolong $(\fg_-,\fg_0)_*$ gives a
Lie algebra that we denote by $11(11;\un;a)$. Its non-positive part
is precisely as that of $\fm\fb(3;\un|8)$ in which we consider the
odd indeterminates even.  Our computation shows that the coordinates
of the shearing vector corresponding to the odd indeterminates can
only be equal to 1, other being arbitrary.

In the $\Zee$-grading considered, the negative part of
$11(11;\un;a)$ is independent of $a$. It first appears in the 0th
component, and only in the three vectors $Z_1$, $H_1=[Z_1, Y_1]$ and
$H_2$. But $Z_1=a\cdot Z'_1$, where $Z'_1$ does not depend on $a$.
So we can replace $Z_1$ by $Z'_1$, and $H_1$ by $[Z'_1, Y_1]$ and
the new basis elements do not depend on $a$.

Finally, replacing $H_2$ by $H_2+H_4$, we get
\[ H_2+H_4=a \,
(x_1\partial_1+   x_3\partial_3 +
   x_4\partial_4 +
   x_5\partial_5 +
   x_7\partial_7 +
   x_9\partial_9 +
   x_{11}\partial_{11})= a\cdot H.
\]
Now, replace $H_2$ by $H$; we get a basis of $\fg_0$ independent of
$a$.

The prolong --- desuperization of $\fm\fb(3;\un|8)$ --- is described
more explicitly in \cite{BGLS}. The realization by vector fields is
as follows:

\footnotesize
$$ \tiny
\renewcommand{\arraystretch}{1.4}
\begin{tabular}{|l|l|} \hline
$\fg_{i}$&the generators  \\ \hline $\fg_{-3}$&$w_{1}=\partial_1,\
w_{2}=\partial_2,\ $\\\hline $\fg_{-2}$&$w_{3}=\partial_3,\
w_{4}=\partial_4,\ w_{5}=\partial_5, $\\\hline
$\fg_{-1}$&$w_{6}=x_5\partial_1+
\partial_6,\  w_{7}=x_5\partial_2+\partial_7,\  w_{8}=x_4\partial_1+
x_7\partial_3+ \partial_8, \  w_{9}=x_4\partial_2+x_6\partial_3+ \partial_9,$\\
&$w_{10}=(x_3+x_6\,
   x_9)\partial_1+x_7\,
   x_9\partial_2+x_7\partial_4+x_9\partial_5+\partial_{10}, \
w_{11}=(x_3+x_7 x_8) \partial_2+ x_6 x_8\partial_1+ x_6
\partial_4+x_8\partial_5 +
\partial_{11}$\\
\hline

$\fg_{0}\simeq $& $Y_1=x_1\partial_2+x_6\partial_7+ x_8\partial_9+ x_{10}\partial_{11},$\\

$\fsl(3)\oplus$&$Z_1= a \, (x_2\partial_1+
   x_7\partial_6 +
   x_9\partial_8 +
   x_{11}\partial_{10}),$\\

$\oplus \fgl(2)$&$ H_1=[Z_1,Y_1]=a \, (x_1\partial_1+  x_2\partial_2
+
   x_6\partial_6+  x_7\partial_7 +
   x_8\partial_8+ x_9\partial_9+
   x_{10}\partial_{10} +  x_{11}\partial_{11})$\\

&$H_3=[Z_3,Y_3]=x_4\partial_4+
   x_5\partial_5  +
   x_6\partial_6  +
   x_7\partial_7  +
   x_8\partial_8 +
   x_9\partial_9 $\\

&$Y_4= x_6\, x_8\, x_9\partial_1+ x_7\, x_8\,
   x_9\partial_2 +
   x_3\partial_4 + x_8\,
   x_9\partial_5 +
   x_8\partial_{10} +
   x_9\partial_{11},$\\

&$ Z_4=x_6\, x_{10}\, x_{11}\partial_1+  x_7\, x_{10}\,
   x_{11}\partial_2 +
   x_4\partial_3+  x_{10}\,
   x_{11}\partial_5 +
   x_{10}\partial_8 +
   x_{11}\partial_9$\\
&$H_4=[Z_4,Y_4]=x_3\partial_3+
   x_4\partial_4 +
   x_8\partial_8 +
   x_9\partial_9 +
   x_{10}\partial_{10} +
  x_{11}\partial_{11},$\\

&$ H_2= a \, x_1\partial_1+ (a
   +1)\, x_3\partial_3 + (a +1)\,
   x_4\partial_4 + a \,
   x_5\partial_5 + a \,
   x_7\partial_7 +
   x_8\partial_8 + (a +1)\,
   x_9\partial_9 +
   x_{10}\partial_{10} + (a +1)\,
   x_{11}\partial_{11}$\\

&$Y_3= x_4\partial_5 +  x_6\,
   x_7 \partial_3 +
 x_6\partial_8 +
 x_7\partial_9 ,\  Z_3= x_5\partial_4 +  x_8\,
   x_9 \partial_3 +
   x_8\partial_6 +
 x_9\partial_7 ,\  Y_7=[Y_3,Y_4],
\ Z_7=[Z_3,Z_4]$
   \\
\hline
\end{tabular}
$$
\normalsize

\sssec{Prolong of the non-positive parts of $\fo_\Pi(2n)/\fc$ and
its derived}

\paragraph{For $r=(10\dots0)$} We consider the Lie
algebras obtained from the algebras with Cartan matrix by
factorizing modulo center and without $\Kee I_0$, i.e.,
$\fo^{(2)}_\Pi(4k+2)$ and $\fo'_\Pi(4k)/\fc$, and in both cases the
0-th part  is isomorphic to $\fo'_\Pi(\dim-2)$, where $\dim= 4k+2$
or $4k$, respectively.

Note that the $4k+2$-dimensional case can be described just like the
$4k$-dimensional one: $\fo^{(2)}_\Pi(4k+2) = \fo'_\Pi(4k+2)/\fc$.
(It is just that it is usually more convenient to use a simpler
description in terms of derived algebra instead of quotient
algebra.) So we can talk about $\fo'_\Pi(\dim)/\fc$ in both cases.

If we consider the 0th part of $\fo'_\Pi(\dim)$ (before
factorization by center) in that grading, it is isomorphic to
$\fo'_\Pi(\dim-2)\oplus \Kee c$, where $c$ is a central element
acting by identity on the $(-1)$st part (take
$d=E^{1,1}+E^{k+1,k+1}$). The algebra $\fo'_\Pi(\dim-2)$ contains
the identity matrix $1_{\dim-2}$ which is also central and also acts
 on the $(-1)$st part by identity. The center of $\fo'_\Pi(\dim)$ modulo which we
factorize the algebra consists of elements that act on the $(-1)$st
part by 0, that is, $\Kee(1_{\dim-2}+ c)$. The resulting quotient
algebra and its action on the $(-1)$st component are the same as if
we factorized modulo $\Kee c$, i.e., it is just $\fo'_\Pi(\dim-2)$.

In table \eqref{rank4} we write $\fo\fc_\Pi(1;8)/\fc$ because we
consider prolongs of (non-positive parts relative a certain
$\Zee$-grading of) Lie algebras with Cartan matrix. In reality we do
not have to first centrally extend an algebra just to factorize it
modulo this center the next moment.

\paragraph{For $2n>8$ and $r=(0\dots01)$}
The prolongation returns $\fo'_\Pi(2n)\ltimes n$ outer derivations.

\paragraph{For $2n>8$ and $r=(0\dots0100)$}
The prolongation returns $\fo'_\Pi(2n)\ltimes (n-2)$ outer
derivations.

\paragraph{For $\fo'_\Pi(8)$, the Dynkin diagram is most symmetric}
For $\fo'_\Pi(8)$, there are
only the two non-equivalent cases:

(a) For the grading $r=(1000)$ (``selected" is any of the end-points
of the Dynkin diagram), the components $\fg_{i}$ for $i\leq 0$ are
as follows:
\begin{equation}\label{endpoint}\tiny
\renewcommand{\arraystretch}{1.4}
\begin{tabular}{|l|l|} \hline
$\fg_{i}$&the generators  \\ \hline \hline $\fg_{-1}$&$\partial_1,\
\partial_2,\ \partial_3,\
\partial_4,\ \partial_5,\ \partial_6$\\
\hline

$\fg_{0}\simeq$& $Y_2=x_1\partial_2+x_5\partial_6;\
Z_2=x_2\partial_1+x_6\partial_5; \ Y_3=x_2\partial_3+
x_4\partial_5; \ Z_3=x_3\partial_2+x_5\partial_4; \ H_5=\mathop{\sum}\limits_{1\leq i\leq 6}x_i\partial_i$;\\
${\bf F}(\fh(0|4))$& $Y_6=[Y_2,Y_3];\ Z_6=[Z_2,Z_3];Y_7=[Y_2,Y_4];\
Z_7=[Z_2,Z_4];\
Y_{10}=[Y_3,Y_7]; \  Z_{10}=[Z_3,Z_7];$\\
$\ltimes \Kee E$ &$H_6=x_3\partial_3+x_5\partial_5+x_6\partial_6;\
H_2=[Z_2,Y_2];\  H_4=[Z_4,Y_4]; \ Y_4=x_2\partial_4+x_3\partial_5;\
Z_4=x_4\partial_2+x_5\partial_3$
   \\
\hline
\end{tabular}
\end{equation} \normalsize
The prolong, denoted by $6(6)$, is of dimension 64 whatever $\un$.
We see that $\dim 6'(6)=62$. The lowest weight vectors of the
$6'(6)_0$-module $6'(6)_1$ are as follows:
\[
\begin{array}{lcl}
v_1 &=&x_1\, x_2\partial_4 + x_1\,
   x_3\partial_5 + x_2\, x_3\partial_6 \\
v_2 &=& x_1\, x_2\partial_+ x_1\,
   x_4\partial_5 + x_2\, x_4\partial_6\\
v_3&=& x_1\, x_2\partial_2 +  x_1\,
   x_5\partial_5 +x_2\, x_5\partial_6 +  x_1\,
   x_3\partial_3 + x_1\, x_4\partial_4+x_3\,
   x_4\partial_6
\end{array}
\]
The modules generated by $v_1$ and $v_2$ are of dimension 7, that
generated by $v_3$ is of dimension 6. Since $\dim\fg_1=20$, these
modules constitute a direct sum (to be sure: there are three highest
weight vectors).

The two vectors that are missing in $6'(6)$ as compared with $6(6)$
act on $6'(6)$ as outer derivatives; one lies in the component
$6(6)_0$, the other one  in  $6(6)_4$; just as in \cite{ILL}.

(b) For the grading $r=(0100)$ (``selected" is the branching node),
we have the following Chevalley basis (the $X_i^\pm$ are Chevalley
generators; of the 4 elements of the maximal torus only $H_1$ and
$H_2$ survive after factorization modulo center):
\tiny
\begin{equation}\label{branch}
\renewcommand{\arraystretch}{1.4}
\begin{tabular}{|l|l|} \hline
$\fg_{i}$&the generators  \\ \hline \hline
$\fg_{-2}$&$w_{1}=\partial_1$\\
\hline $\fg_{-1}$&$X_2^-:=w_2=\partial_2,\ w_3=\partial_3,\
w_4=\partial_4;\ w_5=\partial_5,\ w_6=x_5\partial_1+\partial_6, \
w_7=x_4\partial_1+
\partial_7,$\\
&$ w_{8}=x_3\partial_1+\partial_8;\ w_{9}=x_2\partial_1+\partial_9$\\
\hline $\fg_{0}\simeq$& $X_1^-:=Y_1= x_4\, x_5\partial_1 +
x_2\partial_3 +
   x_4\partial_6 +x_5\partial_7 +x_8\partial_9, \  X_1^+:=Z_1= x_6\,
   x_7\partial_1+ x_3\partial_2 +x_6\partial_4+
   x_7\partial_5 +x_9\partial_8$\\

$\fhei(6)\ltimes$&$X_3^-= Y_3= x_3\, x_5\partial_1 + x_2\partial_4
+x_3\partial_6 +x_5\partial_8 +x_7\partial_9 , \ X_3^+= Z_3= x_6\,
x_8\partial_1 + x_4\partial_2 +x_6\partial_3 +
x_8\partial_5 +x_9\partial_7,$\\

3 outer de-&$X_4^-=Y_4= x_3\, x_4\partial_1+ x_2\partial_5 +x_3\partial_7+
   x_4\partial_8 +x_6\partial_9, \ X_4^+=Z_4= x_7\, x_8\partial_1 + x_5\partial_2+x_7
   \partial_3+
x_8\partial_4 +x_9\partial_6 ,$\\
rivatives&$H_1=[Z_1,Y_1]=[Z_3,Y_3]=[Z_4,Y_4]=\mathop{\sum}\limits_{i\neq
1}x_i\partial_i,\ \ d_3=
H_2=[X_2^+,X_2^-]=\mathop{\sum}\limits_{i\neq 2,6,7,8}x_i\partial_i$
   \\
&$  d_1= x_1\partial_1 +
   x_3\partial_3 +x_6\partial_6 +x_7\partial_7+
   x_9\partial_9, \ d_2= x_1\partial_1 +x_4\partial_4+x_6\partial_6+x_8\partial_8+
x_9\partial_9,$\\
\hline
\end{tabular}
\end{equation} \normalsize
We have $\tilde \fg_{0}:=[\fg_{-1},\fg_1]\simeq\fhei(6)\ltimes
\text{1 outer derivative (which is $H_2$)}$. The $\tilde
\fg_{0}$-module $\fg_1$ is irreducible of dimension 8, with the
lowest weight vector
\[
\begin{array}{lcl}
v:=X_2^+ & =& x_1\, x_2\partial_1 +  x_3\, x_4\,
   x_5\partial_1 +x_1\partial_9 +  x_2\,
   x_3\partial_3 + x_2\, x_4\partial_4+  x_2\,
   x_5\partial_5 +
   x_2\, x_9\partial_9 + \\
&& x_3\,
   x_4\partial_6 +x_3\, x_5\partial_7+  x_3\,
   x_8\partial_9 +
    x_4\, x_5\partial_8+ x_4\,
   x_7\partial_9 +x_5\, x_6\partial_9
\end{array}
\]
Further, set $9(9;\un):=(\fg_{-1},\tilde \fg_{0})_{*,\un}$; by
standard criteria this is a simple Lie algebra; we have
\[
(\fg_{-1},\fg_0)_{*,\un}=9(9;\un)\ltimes\text{ 6 outer derivatives
in the highest component.}
\]
Computer-aided experiments show that $\un=(n,1,\ldots,1)$. The Lie
algebra $9(9;\un)$ is, clearly, an exceptional subalgebra of
$\fk(9;\un)$, a partial prolong.

The Lie algebra $\fo'_\Pi(8)/\fc$ is the result of  desuperization
of an exceptional (in the sense described in eq. (92) in
\cite{BGL1}) simple Lie superalgebra $\fo'_{\Pi\Pi}(4|4)/\fc$. One
can superize $\fo'_\Pi(8)/\fc$ by assuming that any one of the 4
pairs of Chevalley generators is odd; in addition to
$\fo'_{\Pi\Pi}(4|4)/\fc$ this yields $\fo'_{\Pi\Pi}(2|6)/\fc\simeq
\fpe'(4)/\fc$. Equivalently, declaring parities of Chevalley
generators imposes certain restriction on parities of the
indeterminates, see the elements of $\fg_{-1}$ in eq.
\eqref{branch}. So in addition to prolongs of the non-positive part
of $\fo'_{\Pi\Pi}(4|4)/\fc$ which does not differ from that of
$\fo'_{\Pi}(8)/\fc$, except parities of its elements: we declare
$x_3,x_5,x_6,x_8$ odd for $\fo'_{\Pi\Pi}(2|6)/\fc$ and all, except
$x_1$, odd for $\fo'_{\Pi\Pi}(4|4)/\fc$; in this case is
$\fg_0=\fhei(2|4)\ltimes \text{ 3 outer derivatives}$) (or
$\fhei(6)\ltimes \text{ 3 outer derivatives}$, respectively).

The four of the authors suggest to designate this exceptional simple
Lie algebra
\begin{equation}\label{ir}\fii\fr(9;\un):=9(9;\un)\text{~~ and its superizations
$\fii\fr(3;\un|6)$ and $\fii\fr(5;\un|4)$.} \end{equation}

\paragraph{Brown's algebra $D_4(3; \un)$, see \cite{Bro}, as a desuperization of
$\fvle(3;\un|8)$} Set $\cL =\cL(3; \un) =\cL_0\oplus \cL_1$, where
$\cL_0 = \fsvect(3; \un)$ while $\cL_1=\cO((3; \un)_1\oplus \cO(3;
\un)_2$ is the direct sum of two copies of $\cO(3; \un)$ indexed for
convenience, the elements of $\cO(3; \un)_2$ will be barred. Let the
action of $\cL_0$ on $\cL_1$ be the natural one, for any
$f,g\in\cO((3; \un)_1$ and  $\bar f, \bar g\in\cO((3; \un)_2$, set
\[[f,g] = [\bar f,
\bar g] = 0\text{~and $[f, \bar g] = \nabla f \times \nabla g$,}
\]
where $\nabla f = \sum(\partial_if)\partial_i$ and $ D \times  E$ is
determined by bi-linearity over $\cO(3; \un)$ and the rules \[
\begin{array}{ll}
\partial_i\times \partial_i= 0&\text{for $i = 1, 2, 3$}\\
\partial_i\times \partial_j=\partial_k&\text{for $(i, j, k)$ a
permutation of $(1,2,3)$.}\end{array}
\] Brown showed that $\cL$ is a Lie algebra and endowed it with a $\Zee$-grading
as follows: Recall that the standard $\Zee$-grading of the algebra
of functions $\cO(3; \un)$ (degree of each indeterminate equals to
1) induces a $\Zee$-grading (also called standard) of $\fvect(3;
\un)$ and its homogeneous subalgebras, such as $\fsvect(3; \un)$.

Now, let $(\cL_0)_{2i}:=\fsvect(3; \un)_i$ and
$(\cL_1)_{2i-3}:=(\cO((3; \un)_1)_i\oplus (\cO(3; \un)_2)_i$ with
respect to the standard $\Zee$-grading of the algebras in the right
hand sides.

Brown showed that  $\cL_{-3}/\cL_{-2}$ is the center in $\cL$, and
$D_4(3; \un):=\cL/\cL_{-3}$ is a simple Lie algebra; in particular,
we have  $D_4(3; \un_s):=\fo_\Pi^{(2)}(8)/\fc$.

Brown's description of $D_4(3; \un)$ reproduced above can be
formulated in a very simple way as the complete CTS prolong of its
non-positive part, where
\begin{equation}\label{D_4(3}
\fg_0=\fsl(V), \text{~~where $\dim V=3$},\ \fg_{-1}=V\oplus
\overline V, \ \fg_{-2}=V,
\end{equation}
and where $\overline V$ is another copy of the tautological
representation of $\fsl(3)$ whereas the bracket
\begin{equation}\label{crassp}
E^2(\fg_{-1})\tto\fg_{-2}
\end{equation}
is the cross product of vectors of any 3-dimensional space (in other
words, we identify $V$ and $\overline V$ with $\fo'(3)$ and the map
\eqref{crassp} is just the bracket in $\fo'(3)$.

To obtain the Lie algebra $\fsl(3)$ as the 0th part, we have to
consider not the simplest grading, but the one of the form
$r=(0011)$. This grading is not, however, a Weisfeiler one (since
the component $\fg_{-1}$ is not an irreducible $\fg_0$-module) and
the only Weisfeiler regradings of this algebra are the ones from
table \eqref{rank4}. Apart from giving an interpretation of the
Brown's algebra, our answer shows the true number of independent
parameters the vector $\un$ depends on.

Now observe a remarkable likeness of the negative parts of the
following pairs: Brown's $\cL$ and the Lie superalgebra $\fm\fb$,
see subsec. \ref{sssmb}, as well as those of $\cL$ and $\fvle(3|6)$,
see also \cite{ShP}. To prove that $D_4(3; \un)$ is a desuperization
of $\fvle(3|6)$ and $\fvle(3|6)_0$ is contained in $D_4(3; \un)_0$
due to the constraints imposed on $\un$, consider the complete
prolong of the {\bf negative} part of $D_4(3; \un)$ without any
restrictions on $\un$. We see that indeed $D_4(3;
\un)_0=\fsl(3)\oplus\fgl(2)$ which turns $D_4(3; \un)_{-1}$ into an
irreducible $D_4(3; \un)_0$-module and the grading $r=(0011)$ into a
Weisfeiler one. Our computations confirm (Brown's description) that
$\un=(n_1, n_2, n_3, 1,1,1,1,1,1)$. The component $\fg_1$ is of
dimension 12, the sum of 4 irreducible 3-dimensional $D_4(3;
\un)_0$-modules for $\fg=D_4(3; \un)$ with the following lowest
weight vectors (resp. $\dim \fg_1=18$ for $\fg=\fvle(3;\un|6)$, only
the first 2 vectors are generators of the irreducible
$\fvle(3;\un|6)_0$-modules of dimension 8 and 6, respectively):
\[
\begin{array}{ll}
v_1=&x_1\,x_4\partial_2 + x_1\,
   x_6\partial_3+x_4\, x_6\,
   x_7\partial_3+ x_1\partial_9+x_4\,
   x_6\partial_8+x_4\,
   x_7\partial_9\\
v_2=&x_1\,
   x_4\partial_1+ x_2\,
   x_4\partial_2+x_2\,
   x_6\partial_3+ x_1\partial_7+ x_2\partial_9+x_4\,   x_6\partial_6+x_4\,
   x_7\partial_7+x_4\,
   x_8\partial_8+x_4\,
   x_9\partial_9\\
v_3=&x_1\,
   x_5\partial_2+x_1\,
   x_7\partial_3+x_5\, x_6\,
   x_7\partial_3+ x_1\partial_8+ x_5\,
   x_6\partial_8+ x_5\,
   x_7\partial_9\\
v_4=&x_1\,
   x_5\partial_1+ x_2\,
   x_5\partial_2+ x_2\,
   x_7\partial_3+ x_1\partial_6+x_2\partial_8+ x_5\,
   x_6\partial_6+ x_5\,
   x_7\partial_7+ x_5\,
   x_8\partial_8+ x_5\,
   x_9\partial_9
\end{array}
\]
Clearly, Brown's $D_4(3; \un)$ is a partial prolong of the
non-positive part of the above algebra with $\fg_0=\fsl(3)$. The
partial prolong of the non-positive part and just one of the two
irreducible $\fg_0$-modules in $D_4(3; \un_s)_1$ is $\fsvect(3;
\un_s)$. We consider partial prolongs with
$\fg_0=\fsl(3)\oplus\fgl(2)$ in \cite{BGLS}.

We have:
\begin{equation}\label{Broandvle}\tiny
\renewcommand{\arraystretch}{1.4}
\begin{tabular}{|l|l|} \hline
$\fg_{i}$&the generators  \\ \hline \hline
$\fg_{-2}$&$w_{1}=\partial_1,\ w_{2}=\partial_2,\ w_{3}=\partial_3
$\\ \hline $\fg_{-1}$&$w_{4}=\partial_4,\
w_{5}=\partial_5,\  w_{6}=x_5\partial_1+\partial_6,$\\
&$ w_7=x_4\partial_1 +\partial_7,\  w_8=x_5\partial_2+x_7\partial_3+
   \partial_8, \  w_9=x_4\partial_2+x_6\partial_3+
 \partial_9$\\
\hline $D_4(3; \un)_0=\fsl(3)$&$Y_{1}=x_1\partial_2+x_6\,
   x_7\partial_3+ x_6\partial_8+ x_7\partial_9,\  Z_{1}=x_2\partial_1+x_8\,
   x_9\partial_3+ x_8\partial_6+ x_9\partial_7,\ $\\
&$ Y_{2}=x_2\partial_3+x_4\,
   x_5\partial_1+ x_4\partial_6+ x_5\partial_7, \  Z_{2}=x_3\partial_2+x_6\,
   x_7\partial_1+ x_6\partial_4+ x_7\partial_5,\  $\\

&$Y_5=[Y_1, Y_2],\  Z_5=[Z_1, Z_2],\  H_1=[Z_1, Y_1], \
H_2=[Z_2,Y_2] $
   \\
\hline \hline $\fvle(3;\un|6)_{0}=\fsl(3)\oplus \fgl(2)$& the above
with\\
& $\tilde Y_1=x_4\partial_5+
   x_6\partial_7+x_8\partial_9,\, \tilde Z_1=x_5\partial_4+x_7\partial_6+x_9\partial_8, \, \tilde H_1=[\tilde Z_1, \tilde Y_1]$\\
& $\tilde H_2=x_1\partial_1+x_2\partial_2+x_3\partial_3+
 x_5\partial_5+x_7\partial_7+x_9\partial_9$
   \\
\hline
\end{tabular}
\end{equation}\normalsize

\paragraph{$\fvle(4|3;K)$: recapitulation}

If $p=0$, we have $\dim \fg_1=18$, same as above for $p=2$. If
$p=0$, then the $\fg_0$-module $\fg_1$ possesses a 12-dimensional
submodule $V=S^2(\id_{\fsl(3)})\otimes \id_{\fsl(2)}$ (this is what
comes from $\fsle'(3)$, the common parts of the two glued
superalgebras $\fle(3)$ and $\fle(3;3)$, see \cite{ShP}); the
quotient of $\fg_1$ modulo this submodule is of dim 6 and isomorphic
to the tensor product $\id_{\fsl(3)}\otimes \id_{\fsl(2)}$ of
tautological modules over $\fsl(3)$ and $\fsl(2)$. Clearly, for
$p\neq 0$ (except, perhaps, for $p=3$; we have to check), both $V$
and the quotient are irreducible.  The quotient module is not a
direct summand: $\fg_1$ is an indecomposable $\fg_0$-module.

Let us denote the indeterminates of functions that generate the two
copies of $\fle(3)$ with common part $\fsle'(3)$ by $u_1, u_2,u_3,
\xi_1,\xi_2,\xi_3$ and $u'_1, u'_2,u'_3, \xi'_1,\xi'_2,\xi'_3$,
respectively, assuming that if
\[
\begin{array}{l}
\deg_\xi f(u,\xi)=1\text{~~and~~}\Delta f=0,\\
\text{harmonic functions being singled out by the ``odd
Laplacian"~~}\Delta =\sum
\frac{\partial^2}{\partial_{u_i}\partial_{\xi_i}};
\end{array}
\]
then we identify
\[
\begin{array}{l}
f(u,\xi)\text{~~with~~}f(u',\xi'),\\
f(u)\text{~~with~~}\sum \frac{\partial
f(u')}{\partial_{u'_i}}\xi'_j\xi'_k \text{~~for any even permutation
$(i,j,k)$ of $(1,2,3)$}.\end{array}
\]
Now, set $\deg\xi_i=\deg\xi'_i=1$ and $\deg u_i=\deg u'_i=2$ for all
indices, whereas the degree of the element of $\fvle$ with
generating function $f$ is $\deg_{Lie}(f)=\deg(f)-3$. Since we
factorize both copies of the spaces of generation functions modulo
constants, the depth of the resulting Lie superalgebra is equal to
2.

The component $\fg_{-2}$ is spanned by functions of degree 1, i.e.,
by $\xi_i=\xi'_i$ (3 elements).

The component $\fg_{-1}$ is spanned by functions of degree 2, i.e.,
$u_i=\xi'_i\xi'_k$ and $\xi_i\xi_j=u'_k$ for any even permutation
$(i,j,k)$ of $(1,2,3)$.

The component $\fg_0$ is spanned by functions of degree 3. These are
\[\fsl(3)=\{f\in \Span(u_i\xi_j)\mid \Delta f=0\}=\{f\in
\Span(u'_i\xi'_j)\mid \Delta f=0\},
\]
and also
\[
\fgl(2)=\Span(\sum u_i\xi_i,\ \ \sum u'_i\xi'_i,\ \
\xi_1\xi_2\xi_3,\ \ \xi'_1\xi'_2\xi'_3).
\]

The component $\fg_1$ is spanned by degree 4 functions: 6 monomials
of degree 2 in $u$ span $V$, and 9 monomials of the form
$u_i\xi_\alpha\xi_\beta$ span a subspace $W$; analogous monomials in
primed indeterminates span $V'$ and $W'$. How to glue these spaces?

We have $W\supset W_0$ (and $W'\supset W'_0$), where subspaces $W_0$
and $W'_0$ consist of harmonic functions. We see (from \cite{ShP})
that $V$ is glued with $W'_0$, while $W_0$ with $V'$. The subspace
$V\oplus W_0=V'\oplus W'_0$ obtained is precisely
$S^2(\id_{\fsl(3)})\otimes \id_{\fsl(2)}$ as $\fg_0$-module; it is
exactly component of degree 1 of the subalgebra $\fsle'$ in the
grading considered.

The quotient of $\fg_1$ modulo this submodule is a $\fg_0$-module of
dimension 6 isomorphic to $\id_{\fsl(3)}\otimes \id_{\fsl(2)}$; the
$\fg_0$-module $\fg_1$ is indecomposable.

\paragraph{About $\fvle(4|3;1)$} Let $\fg$ be a Lie algebra,
$\cO(x)$ an associative algebra of ``functions" in indeterminates
$x$ (polynomials, divided powers, etc.), and $d$ a derivation of
$\fg$. The expression
$$
X\mapsto 1\otimes X+d\otimes \Div(X), \text{~~where $X\in
\fvect(x)$},
$$
determines a $\fvect(x)$-action on $\fh=\fg\otimes \mathcal O(x)$
commuting with the operator $d\otimes 1$.

Now, if we identify the $(-1)$st component of  $\fvle(4|3;1)$ with
$V\otimes \Lambda(2)$, then the 0th component would contain an ideal
$\fsl(2)\otimes\Lambda(2)$, a subalgebra isomorphic to $\fvect(2)$,
acting on the  $(-1)$st component as $1\otimes X+E_{11}\otimes
\Div(X)$, where $X\in\fvect(2)$, and commuting with the ideal, and
instead of the center we have to add the derivation $E_{11}\otimes
1$, where $E_{11}$ is a matrix unit.

\begin{landscape}
\begin{table}[ht]\centering
{\footnotesize%
\caption{Dynkin diagrams for Lie superalgebras: $p=2$}
\begin{equation}\label{tbl}
\extrarowheight=2pt \begin{tabular}%
{|>{\PBS\raggedright\hspace{0pt}}m{33mm}|
 >{\PBS\raggedright\hspace{0pt}}m{40mm}|
 >{\PBS\raggedright\hspace{0pt}}m{13mm}|
 >{\PBS\raggedright\hspace{0pt}}m{17mm}|
 >{\PBS\raggedright\hspace{0pt}}m{13mm}|
 >{\PBS\centering\hspace{0pt}}m{12mm}|
 >{\PBS\raggedright\hspace{0pt}}m{25mm}|}
\hline
 Diagrams\centering & $\fg$\centering & $v$\centering & $ev$ &
 $od$ &$png$ & $ng\le \min(* \ , *)$\\ \hline
 &&&
 $k_\ev-2$ & $k_\od$ &$\ev$ &$2k_\ev-4, 2k_\od$ \\
 &&&
 $k_\od$ &$k_\ev-2$ & $\od$ &$2k_\ev-3, 2k_\od-1$ \\
 &&&
 $k_\od-2$ & $k_\ev$ &$\ev$ &$2k_\ev, 2k_\od-4$ \\
 &&&
 $k_\ev$ &$k_\od-2$ & $\od$ &$2k_\ev-1,2k_\od-3$ \\ \cline{4-7}
 &&&
 $k_\ev-1$ & $k_\od-1$ & &$2k_\ev-2, 2k_\od-1$ \\
 \raisebox{3.5em}[0pt]{$\left.\arraycolsep=0pt\begin{array}{l}
 1)\ \includegraphics{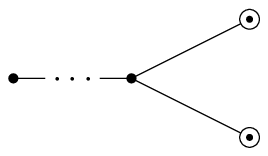}\\
 2)\ \includegraphics{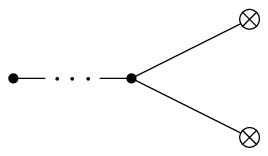}
 \end{array}\ \right\}$} &
 \raisebox{3.5em}[0pt]{$\begin{array}{l}\foo\fc(2;2k_\ev|2k_\od)\ltimes\Kee I_0\\
 \text{if $k_\ev+k_\od$ is odd;}\\ \foo\fc(1;2k_\ev|2k_\od)\ltimes\Kee I_0\\
 \text{if $k_\ev+k_\od$ is even.}\end{array}$}
 &
 \raisebox{3.5em}[0pt]{$k_{\ev}+k_{\od}$} &
 $k_\od-1$ &$k_\ev-1$ & &$2k_\ev-1, 2k_\od-2$ \\
\hline
 &&&
 $k_\ev-1$ & $k_\od$ &$\ev$ &$2k_\ev-2, 2k_\od$ \\
 &&&
 $k_\od$ &$k_\ev-1$ & $\od$ &$2k_\ev-1,2k_\od-1$ \\ \cline{4-7}
 &&&
 $k_\od-1$ & $k_\ev$ & $\ev$ &$2k_\ev, 2k_\od-2$ \\
 \raisebox{2.2em}[0pt]{$\left.\arraycolsep=0pt\begin{array}{l}
 3)\  \includegraphics{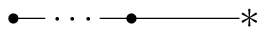}\\
 4)\ \includegraphics{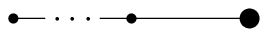}
 \end{array}\ \right\}$} &
 \raisebox{2.2em}[0pt]{$\foo'_{I\Pi}(2k_\ev+1|2k_\od)$} &
 \raisebox{2.2em}[0pt]{$k_{\ev}+k_{\od}$} &
 $k_\ev$ &$k_\od-1$ & $\od$ &$2k_\ev-1, 2k_\od-1$ \\
 \hline
 5)\ \includegraphics{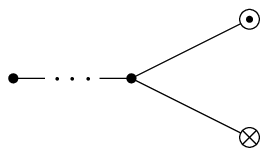}
 & $\begin{array}{l}\fpe\fc(2;m)\ltimes\Kee I_0\\ \text{for $m$ odd;}\\
 \fpe\fc(1;m)\ltimes\Kee I_0\\ \text{for $m$ even.}\end{array}$ & $m$&&&&\\ \hline
\end{tabular}
\end{equation}
\vskip 1 cm {\bf Notation} The Dynkin diagrams in Table \ref{tbll}
correspond to Cartan matrix Lie superalgebras close to
ortho-orthogonal and periplectic Lie superalgebras. Each thin black
dot may be $\motimes$ or $\odot$; the last five columns show
conditions on the diagrams; in the last four columns, it suffices to
satisfy conditions in any one row. Horizontal lines in the last four
columns separate the cases corresponding to different Dynkin
diagrams. The notation are:
 $v$ is the total number of nodes in the diagram;  $ng$ is the
 number of ``grey'' nodes $\motimes$'s among the
thin black dots;  $png$ is the parity of this number;
 $ev$ and $od$ are the number of thin black dots such that the number
of $\otimes$'s to the left from them is even and odd, respectively.
\vspace{-4mm} }
\end{table}
\end{landscape}

\end{document}